\documentclass[12pt,a4paper]{amsart}
\usepackage{a4}
\usepackage[english]{babel} 
\usepackage{amsmath,amssymb,amsthm}
\def\leq {\leqslant}
\def\le {\leqslant}
\def\ge {\geqslant}
\def\geq {\geqslant}
\def\@bibitem[#1]#2{\item\@biblabel{#1}.\if@filesw
{\def\protect##1{\string##1\space}\immediate\write
\@auxout{\string\bibcite{#2}{#1}}}\fi\ignorespaces\@showtag{#2}}

\theoremstyle{plain}
\newtheorem{theorem}{Theorem}
\newtheorem{rem}{Remark}
\newtheorem{lemma}{Lemma}
\renewcommand{\theequation}%
{\arabic{section}.\arabic{equation}}
\begin{document}

\title{ On exact estimates of the order of approximation of functions of several variables in the anisotropic Lorentz - Zygmund space}
\author{ G. Akishev}
\address{ Lomonosov Moscow University, Kazakhstan Branch \\
Str. Kazhymukan, 11 \\
010010, Nur-Sultan, Kazakhstan}

\address{
Institute of mathematics and mathematical modeling\\
Pushkin str, 125 \\
050010, Almaty \\
 Republic of Kazakhstan
 }
 
\address{Ural Federal University,
Institute of Natural Sciences and Mathematics,
 Yekaterinburg,  Russia 
}

\maketitle

\begin{quote}
\noindent{\bf Abstract.}
In this paper we consider $L_{\overline{p}, \overline\alpha, \overline{\tau}}^{*}(\mathbb{T}^{m})$ anisotropic Lorentz-Zyg\-mu\-nd space $ 2\pi$ of periodic functions of $m$ variables and Nikol'skii--Besov's class $S_{\overline{p}, \overline\alpha, \overline{\tau}, \bar{\theta}}^{\bar r}B$.  
In this paper, we establish order-sharp estimates of the best approximation by trigonometric polynomials with harmonic numbers from the step hyperbolic cross of functions from the Nikol'skii - Besov class in the norm of the anisotropic Lorentz-Zygmund space.
\end{quote}
\vspace*{0.2 cm}

{\bf Keywords:} Lorentz-Zygmund space, \and  Nikol'skii-Besov class, \and trigonometric polynomial,\and best  approximation

 {\bf MSC:} 41A10 and 41A25,  42A05


\section*{Introduction}


Let $\mathbb{N}$, $\mathbb{Z}$, $\mathbb{R}$ be sets of natural, integer, real numbers, respectively, and $\mathbb{Z}_{+} = \mathbb{N}\cup\{0\}$,
$\mathbb{R}^{m}$ --- $m$ -- dimensional Euclidean point space $\bar{x} = (x_{1}, \ldots, x_{m})$ with real coordinates; $I^{m} = \{\bar{x} \in \mathbb{R}^{m}; \ 0 \leq x_{j} \leq 1; \ j = 1, \ldots, m \}=[0, 1]^{m} $ --- $m$--dimensional
cub and $\mathbb{Z}_{+}^{m}$ is the Cartesian product of the sets $\mathbb{Z}_{+}$, $m$ is times.

 \smallskip

Let numbers $p, \tau \in (1, \infty),$
$\alpha \in \mathbb{R}.$ Lorentz-Zygmund space  $L_{p, \tau}\left(\log L \right)^{\alpha}(\mathbb{T})$ is the set of all Lebesgue measurable and $2 \pi$ periodic functions $ f $ for which (see e.g. \cite{1})
$$
\|f\|_{p, \alpha, \tau } := \left\{\int_{0}^{1}\biggl(f^{*}
(t)\biggr)^{\tau}\biggl(1 + |\log_{2} t|\biggr)^{\alpha \tau}t^{\frac{\tau}{p}-1}dt
\right\}^{\frac{1}{\tau}} < +\infty.
$$
where $f^{*}(t)$ is a non-increasing rearrangement of the function
$|f(2\pi x)|$, $x \in [0, 1]$, $\mathbb{T}=[0, 2\pi]$.

The Lorentz-Zygmund space is sometimes denoted by the symbol $L_{p, \alpha, \tau }(\mathbb{T})$. We will use this notation.

Note that for $\alpha=0$ the space $L_{p, \alpha, \tau }(\mathbb{T})$ coincides with the Lorentz space $L_{p, \tau}(\mathbb{T})$, $1 <p, \tau <\infty$, which consists of all functions $f$ such that (see e.g. \cite[Ch. 1, Sec. 3]{2})
$$
\|f\|_{p, \tau}^{*}:=\left(\int\limits_{0}^{1}f^{*^{\tau}}(t)t^{\frac{\tau}{p}-1}dt\right)^{1/\tau} < \infty.
$$

It is known that the fundamental function $\varphi_{L_{p, \tau, \alpha}}(t)$ of the Lorentz-Zygmund  $L_{p, \alpha, \tau}(\mathbb{T})$ satisfies the relation $\varphi_{L_{p, \tau, \alpha}}(t) \asymp
t^{\frac{1}{p}}(1 + |\log_{2} t|)^{\alpha}, t\in (0, 1]$.

Let 
 $\bar{p}=(p_{1},\ldots p_{m})$, $\bar{\tau}=(\tau_{1},\ldots \tau_{m}),$ $\bar{\alpha}=(\alpha_{1},\ldots \alpha_{m})$ and $p_{j}, \tau_{j} \in (1, \infty),$ $\alpha_{j} \in \mathbb{R}$, $j=1,...,m$.
By $L_{\overline{p}, \overline\alpha, \overline{\tau}}^{*}(\mathbb{T}^{m})$ we denote the anisotropic Lorentz-Zygmund space of all measurable Lebesgue functions of $m$ variables $f$ with period
$2 \pi$ for each variable and for which the quantity 
$$
\|f\|_{\bar{p}, \bar{\alpha}, \bar{\tau}}^{*} := \|\ldots\|f^{*_{1},...,*_{m}}\|_{p_{1}, \alpha_{1}, \tau_{1}}\ldots\|_{p_{m}, \alpha_{m}, \tau_{m}}
$$
$$
\Bigl[\int_{0}^{1}
\Bigl[\ldots\Bigl[\int_{0}^{1}\left(
f^{*_{1},...,*_{m}}(t_{1},...,t_{m})
\right)^{\tau_{1}}\biggl(\prod_{j=1}^{m}\biggl(1 + |\log_{2} t_{j}|\biggr)^{\alpha_{j}}t_{j}^{\frac{1}{p_{j}}-\frac{1}{\tau_{j}}}\biggr)^{\tau_{1}}dt_{1}\Bigr]^
{\frac{\tau_{2}}{\tau_{1}}} \ldots
\Bigr]^{\frac{\tau_{m}}{\tau_{m-1}}}dt_{m}
\Bigr]^{\frac{1}{\tau_{m}}}
$$
is finite, where
$f^{*_{1},...,*_{m}}(t_{1},...,t_{m})$
non-increasing rearrangement  of a function
$|f(2\pi \bar{x})|$ for each variable $x_{j} \in [0, 1]$ with fixed other variables (see \cite{3}, \cite{4}). Here and in what follows, $\mathbb{T}^{m} = [0, 2\pi]^{m}$.

For $\alpha_{j}=0$, $j=1,...,m$  space $L_{\overline{p}, \overline\alpha, \overline{\tau}}^{*}(\mathbb{T}^{m})$ is an anisotropic Lorentz space and is denoted by $L_{\overline{p}, \overline{\tau}}^{*}(\mathbb{T}^{m})$
, and $\|f\|_{\bar{p}, \bar{\alpha}, \bar{\tau}}^{*} =\|f\|_{\bar{p},  \bar{\tau}}^{*}$ (see \cite{5}).

If $\alpha_{j}=0$ and $p_{j} = \tau_{j}=p$ , $j=1,...,m$,
then $L_{\overline{p}, \overline\alpha, \overline{\tau}}^{*}(\mathbb{T}^{m})= L_{p}(\mathbb{T}^{m})$ is the well-known Lebesgue space with the norm
$$
\|f\|_{q}:=\biggl(\,\int\limits_{I^{m}}|f(2\pi\bar{x})|^{q}d\bar{
x}\,\biggr)^{1/q},\quad 1\le q < \infty.
$$

For a sequence of numbers we write $\left\{ a_{\overline{n} } \right\}
_{\overline{n} \in \mathbb{Z}_{0}^{m} } \in l_{\overline{p} }$ if

\begin{equation*}  
\Bigl\|
\Bigl\{a_{\overline{n}}\Bigr\}_{\bar{n} \in \mathbb{Z}_{+}^{m}}
\Bigr\|_{l_{\overline{p}}(\mathbb{Z}_{+}^{m})} = \Bigl\{\sum\limits_{n_{m}
=0}^{\infty }\Bigl[ ...\Bigl[\sum\limits_{n_{1} =0
}^{\infty }\Bigl|a_{\overline{n}}\Bigr|^{p_{1}}\Bigr]
^{\frac{p_{2}}{p_{1}}} ...\Bigr]^{\frac{p_{m}}{p_{m-1}}}
\Bigr\}^{\frac{1}{p_{m}}} <+\infty ,
\end{equation*} 
where $\overline{p} =\left( p_{1} ,...,p_{m} \right)$, $1\leq
p_{j} <+\infty$, $j=1,2,...,m$.

If $p_{j}=\infty, \;\; j=1,...,m$, then  
$$
\Bigl\|\{a_{\bar{n}}\}\Bigr\|_{l_{\infty}(\mathbb{Z}_{+}^{m})}=\sup\limits_{\bar{n}\in\mathbb{Z}_{0}^{m}}|a_{\bar{n}}|.
$$

We will use the following notation:
let $\overset{\circ \;\;}
L_{\bar{p}, \bar\alpha, \bar{\tau}}^{*}
\left(\mathbb{T}^{m} \right)$ be the set of functions $f\in
L_{\bar{p}, \bar\alpha, \bar{\tau}}^{*}(\mathbb{T}^{m})$ such
that
 $$
\int\limits_{0}^{2\pi }f\left(\overline{x} \right) dx_{j}
=0,\;\;\forall j=1,...,m;
 $$
$a_{\overline{n} } (f)$ be the Fourier coefficients of $f\in
L_{1}(\mathbb{T}^{m})$ with respect to the multiple
trigonometric system and 
$$
\delta _{\overline{s} } \left( f,\overline{x} \right)
:=\sum\limits_{\overline{n} \in \rho \left( \overline{s} \right)
}a_{\overline{n} } \left( f\right) e^{i\langle\overline{n} ,\overline{x}\rangle } ,
$$
where $\langle\bar{y}, \bar{x}\rangle=\sum\limits_{j=1}^{m}y_{j}
x_{j}$,
 $$
\rho (\bar{s}):=\left\{ \overline{k} =\left( k_{1}
,...,k_{m} \right) \in \mathbb{Z}^{m}: \quad 2^{s_{j} -1} \leq \left| k_{j}
\right| <2^{s_{j} } ,j=1,...,m\right\}.
$$ 
 
  Let $\bar{\gamma}=(\gamma_{1},\ldots,\gamma_{m}),$ $\gamma_{j}>0, j=1,\ldots,m$.
$E_{n}^{(\overline{\gamma})}(f)_{\bar{p}, \bar{\alpha}, \bar{\tau}}$ is the best approximation of the function $f\in L_{\overline{p}, \overline\alpha, \overline{\tau}}^{*}(\mathbb{T}^{m})$ trigonometric polynomials with harmonic numbers from the set $Q_{n}^{\bar\gamma}=
 \cup_{{}_{\langle\bar{s},\bar{\gamma}\rangle < n}}\rho(\bar{s})$ is a stepped hyperbolic cross. For $\bar{\alpha}=(0, \ldots , 0)$ and $\bar{\tau}=\bar{p}$ the  quantity is $E_{n}^{(\overline{\gamma})}(f)_{\bar{p}, \bar{\alpha}, \bar{\tau}}$ is denoted as $E_{n}^{(\overline{\gamma})}(f)_{\bar{p}}$ (see \cite{12} and references therein).
 
In function theory, the $S_{p, \theta}^{\bar r}B$ space of Nikol'skii-Besov in the Lebesgue space $L_{p}$, $1\leq p < \infty$ and its applications is well known (see  \cite{5}-\cite{8}).

In this article, we will consider an analogue of the Nikol'skii-Besov class in the anisotropic Lorentz-Zygmund space:
 $$
S_{\bar{p}, \bar\alpha, \bar{\tau}, \bar{\theta}}^{\bar r}B :=
\Bigl\{f\in \overset{\circ \;\;}
L_{\bar{p}, \bar\alpha, \bar{\tau}}^{*}
\left(\mathbb{T}^{m} \right) : \,\,
\|f\|_{\bar{p}, \bar{\alpha}, \bar{\tau}}^{*} + \Bigl\|\Bigl\{\prod_{j=1}^{m}
2^{s_{j}r_{j}} \|\delta_{\bar
s}(f)\|_{\bar{p}, \bar{\alpha}, \bar{\tau}}^{*} \Bigr\}_{\bar{s}\in
\mathbb{Z}_{+}^{m}}\Bigr\|_{l_{\bar\theta}}\leq 1\Bigr\},
 $$
where $\bar{p}=(p_{1},...,p_{m})$, $\bar\alpha=(\alpha_{1},...,\alpha_{m})$,  $\bar{\tau}=(\tau_{1},...,\tau_{m})$, $\bar{\theta}=(\theta_{1},...,\theta_{m}),$  $\bar{r}=(r_{1},...,r_{m}),$ $1<p_{j}, \tau_{j}< \infty$, 
$0 < \theta_{j}\le+\infty,$ $0<r_{j}<+\infty,$ $\alpha_{j}\in \mathbb{R}$, $j=1,...,m.$

We put  
$$
E_{n}^{(\overline{\gamma})}(S_{\bar{p}, \bar\alpha, \bar{\tau}^{(1)},  \bar{\theta}}^{\bar r}B)_{\bar{q}, \bar{\beta}, \bar{\tau}^{(2)}} := \sup\limits_{f \in S_{\bar{p}, \bar\alpha, \bar{\tau}^{(1)}, \bar{\theta}}^{\bar r}B}E_{n}^{(\overline{\gamma})}(f)_{\bar{q}, \bar{\beta}, \bar{\tau}^{(2)}},
$$
where $\bar{q}=(q_{1},\ldots , q_{m})$, $\bar{\beta}=(\beta_{1},\ldots , \beta_{m})$, $\bar{\tau}^{(i)}=(\tau_{1}^{(i)},\ldots , \tau_{m}^{(i)})$ and $1<q_{j}, \tau_{j}^{(i)}< \infty$, $i=1, 2$, $\beta_{j}\in \mathbb{R}$, $j=1,...,m$.   

Order-exact estimates $E_{n}^{(\overline{\gamma})}(f)_{\bar{p}}$ - the best approximation of functions from the Sobolev classes $W_{p}^{\bar r}$ and Nikol'skii - Besov $S_{p, \theta}^{\bar r}B$, in the metric of the space $L_{q}(\mathbb{T}^{m})$, $1< p, q < \infty$ are well known and are given in the survey articles \cite{9}, \cite{10} and in the monographs \cite{11}, \cite{12} (for more details, see the references therein).
These questions in the anisotropic Lorentz space $L_{\bar{q},\bar{\tau}}^{*}(\mathbb{T}^{m})$ were studied in \cite{13}-\cite{18}.

The main aim of the present paper is to find the order of the quantity 
$$
E_{n}^{(\overline{\gamma})}(S_{\bar{p}, \bar\alpha, \bar{\tau}^{(1)},  \bar{\theta}}^{\bar r}B)_{\bar{q}, \bar{\beta}, \bar{\tau}^{(2)}}.
$$

This paper is organized as follows. In  Section 1 we give auxiliary results.  
 In Section 2, we present and prove the main results.  
 
 We shall denote by  $C(p,q,y,..)$  positive quantities which depend
only on the parameter in the parentheses and not necessarily the
same in distinct formulas .
 The notation $A\left( y\right) \asymp
B\left( y\right)$ means that there exist positive constants
 $C_{1},\,C_{2} $ such that  $C_{1} A(y)
\leq B(y) \leq C_{2} A(y)$.
 For brevity, in the case of the inequalities $ B \geq C_{1} A$ or $B \leq C_{2}A$, we often write $B >> A$ or $B << A $, respectively.

\section{Auxilary statements}\label{sec 1}

In this section, we prove several lemmas necessary to prove the main results of the article.

\begin{lemma}\label{lem 1} 
Let $ \alpha, \beta \in \mathbb{R}$, $l\in \mathbb{N}$.  

1) If $1- \alpha > 0$ and $1-\beta >0$, then 
$$ 
\sum\limits_{0 \leq s < l} (s+1)^{-\alpha}(l-s)^{-\beta} \asymp (l+1)^{1-(\alpha+\beta)}. 
$$

2) If $\alpha = 1, \beta = 1,$ then
 $$
\sum\limits_{0 \leq s < l} (s+1)^{-1}(l-s)^{-1} \asymp l^{-1}\ln (1+l). 
 $$

3) If $\alpha = 1, 0 < \beta <  1,$ then 
 $$
\sum\limits_{0 < s < l} s^{-1}(l-s)^{-1} \asymp l^{-\beta}\ln (1+l). 
 $$
\end{lemma}
{\bf Proof.}  If $ \alpha, \beta \in \mathbb{R}$, then 
$$
\sum\limits_{0 \leq s < l} (s+1)^{-\alpha}(l-s)^{-\beta} =  \sum\limits_{0 \leq s < \frac{l}{2}} (s+1)^{-\alpha}(l-s)^{-\beta} +
 \sum\limits_{\frac{l}{2} \le s < l}(s+1)^{-\alpha}
(l-s)^{-\beta}.                \eqno (1) 
$$
If $1-\alpha > 0,$ then 
$$
C_{1}(\alpha)(l+1)^{1-\alpha}\le \sum\limits_{0 \leq s < \frac{l}{2}} (s+1)^{-\alpha}
\le C_{2}(\alpha)(l+1)^{1-\alpha}.       \eqno (2)  
$$
If $1-\beta > 0,$ then
$$
C_{1}(\beta)(l+1)^{1-\beta} \le \sum\limits_{\frac{l}{2} \le s < l}
(l-s)^{-\beta} \le C_{2}(\beta)(l+1)^{1-\beta}.   \eqno (3)  
$$

If $\beta < 0$ and $1-\alpha > 0,$ then
$$
\sum\limits_{0 \leq s < l} s^{-\alpha}(l-s)^{-\beta} \le l^{-\beta} \sum\limits_{0 \leq s < l} s^{-\alpha}\leq C_{2}(\alpha)(l+1)^{1-\alpha-\beta}.  \eqno (4)
$$

If $\alpha < 0$ and  $1 - \beta  > 0,$ then
$$
\sum\limits_{0 \leq s < l} s^{-\alpha}(l-s)^{-\beta} \le l^{-\alpha} \sum\limits_{0 \leq s < l} (l-s)^{-\beta}\leq C(\alpha)(l+1)^{1-\alpha-\beta}.  \eqno (5)
$$
From inequalities (1)--(5), we obtain
 
$$
\sum\limits_{0 \leq s < l} (s+1)^{-\alpha}(l-s)^{-\beta} \asymp 
(l+1)^{1-(\alpha+\beta)}.
$$
The first point is proven.

2) If $\alpha=\beta=1,$ then 
$$
\sum\limits_{0 \leq s < \frac{l}{2}} (s+1)^{-1} \asymp \ln(1+l)  \eqno (6)
$$
and
$$
\sum\limits_{\frac{l}{2} \le s < l}
(l-s)^{-1} \asymp \ln(1+l).  \eqno (7)  
$$
Therefore, according to inequalities (6), (7), we have
$$
\sum\limits_{0 \leq s < l} (s+1)^{-1}(l-s)^{-1}\le \left(\frac{l}{2}\right)^{-1}
\sum\limits_{0 \leq s < \frac{l}{2}} (s+1)^{-1} + \left(\frac{l}{2}+1\right)^{-1}
\sum\limits_{\frac{l}{2} \le s < l}(l-s)^{-1} \le C_{3}l^{-1}\ln (1+l).
$$
Reverse inequality
$$
\sum\limits_{0 \leq s < l} (s+1)^{-1}(l-s)^{-1}\ge \sum\limits_{0 < s < l} (s+1)^{-1}l^{-1}\ge Cl^{-1}\ln (1+l).
$$
The second statement is proved.

3) If $\alpha = 1, 0 <\beta <1, $ then from (5) we obtain 
$$
\sum\limits_{0 < s < l} s^{-1}(l-s)^{-\beta} \le \left(\frac{l}{2}\right)^{-\beta} 
\sum\limits_{0 < s < \frac{l}{2}} s^{-1} +
\left(\frac{l}{2}\right)^{-1} \sum\limits_{\frac{l}{2} \le s < l}(l-s)^{-\beta}
\le
$$
$$
\le Cl^{-\beta}\ln(1+l).
$$
Lemma 1 is proved

\begin{lemma}\label{lem 2} 
Let $\theta, \beta \in (0, \infty)$. Then, for any  $\lambda_{1}, \lambda_{2} \in \mathbb{R}$, the following relation holds: 
$$
\sum\limits_{s=0}^{n} 2^{-s\beta\theta}(s+1)^{\lambda_{2}\theta}(n-s+1)^{\lambda_{1}\theta} \asymp (n+1)^{\lambda_{1}\theta} \eqno (8) 
$$
$$
\sum\limits_{s=0}^{n} 2^{s\beta\theta}(s+1)^{\lambda_{2}\theta}(n-s+1)^{\lambda_{1}\theta} \asymp 2^{n\beta\theta}(n+1)^{\lambda_{2}\theta}  
\eqno (9)
$$
 \end{lemma} 
{\bf Proof.} Let $n\in \mathbb{N}$, $n\leq x < n+1$. Then
$$
F(x) = \int_{0}^{x}2^{-t\beta\theta}(t+1)^{\lambda_{2}\theta}(x-t+1)^{\lambda_{1}\theta}dt\geq 
$$
$$
\sum\limits_{k=0}^{n-1} 2^{-k\beta\theta}(k+1)^{\lambda_{2}\theta}(n-k+1)^{\lambda_{1}\theta}. \eqno (10)
$$
Therefore, it suffices to estimate the function $F(x)$.

Let 
$\beta > 0$. We will choose a number 
  $\delta\in (0, \beta)$,   
then 
$$
F(x) = 2^{-x\theta\delta}\int_{0}^{x}2^{-t(\beta-\delta)\theta}(t+1)^{\lambda_{2}\theta} 2^{(x-t)\theta\delta}(x-t+1)^{\lambda_{1}\theta}dt.
$$
Since $0\leq x-t\leq x$ and the function $2^{y\delta}(y+1)^{\lambda_{1}}$ increases by $ [0, \infty)$, then
$$
F(x)\leq 2^{-x\theta\delta} 2^{x\theta\delta_{1}}(x+1)^{\lambda_{1}\theta}\int_{0}^{x}2^{-t(\beta-\delta)\theta}(t+1)^{\lambda_{2}\theta}dt
$$
$$
\leq (x+1)^{\lambda_{1}\theta}\int_{0}^{\infty}2^{-t(\beta-\delta)\theta}(t+1)^{\lambda_{2}\theta}dt. \eqno (11)
$$
Since $\beta-\delta> 0$, the integral 
$$
\int_{0}^{\infty}2^{-t(\beta-\delta)\theta}(t+1)^{\lambda_{2}\theta}dt
$$
converges. 
Hence from (11), we obtain
$$
F(x)\leq C(\beta, \theta, \lambda_{2})(x+1)^{\lambda_{1}\theta}
$$
for  $\lambda_{1}, \lambda_{2}\in \mathbb{R}$.

Let us now prove that this estimate is sharp.
If $\lambda_{1}< 0$, then $(x-t+1)^{\lambda_{1}}\geq (x+1)^{\lambda_{1}}$ for $0\leq t \leq x$.
Therefore, 
$$
F(x)\geq (x+1)^{\lambda_{1}\theta}\int_{0}^{1/2}2^{-t\delta_{1}\theta}(t+1)^{\lambda_{2}\theta}dt
$$
for  $x\geq 1/2$.

If $\lambda_{1}> 0$, then
$$
F(x)\geq \int_{0}^{x/2}2^{-t\beta\theta}(t+1)^{\lambda_{2}\theta}(x-t+1)^{\lambda_{1}\theta}dt\geq (\frac{x}{2}+1)^{\lambda_{1}\theta}\int_{0}^{x/2}2^{-t\beta\theta}(t+1)^{\lambda_{2}\theta}dt
$$
$$
\geq (\frac{x}{2}+1)^{\lambda_{1}\theta}\int_{0}^{1/4}2^{-t\beta\theta}(t+1)^{\lambda_{2}\theta}dt
$$
for $x\geq 1/2$.

Now relation (8) follows from these inequalities and (10), (11).

Let $\beta < 0$. We will choose the number $\eta \in (0, -\beta)$. 
Now taking into account that the function $2^{t\eta}(t+1)^{\lambda_{2}}$ is increasing and making the change of the variable $y = x-t$, we get
 $$
F(x) \leq 2^{x\theta\eta}(x+1)^{\lambda_{2}\theta}\int_{0}^{x}2^{t(-\beta-\eta)\theta}(x-t+1)^{\lambda_{1}\theta}dt=
$$
$$
= 2^{x\theta\eta}(x+1)^{\lambda_{2}\theta}2^{x(-\beta-\eta)\theta}\int_{0}^{x}2^{-y(-\beta-\eta)\theta}(y+1)^{\lambda_{1}\theta}dy
$$
$$
\leq (x+1)^{\lambda_{2}\theta}2^{-x\beta\theta}\int_{0}^{\infty}2^{-y(-\beta-\eta)\theta}(y+1)^{\lambda_{1}\theta}dy.
$$
It is not difficult to verify the accuracy of this estimate. The lemma is proved.

Next, we consider the sets
$Y^{m}(\bar{\gamma}, n)=\{\bar{s}\in\mathbb{Z}_{+}^{m}: \,\,
 \langle\bar{s},\bar{\gamma}\rangle\geq n \}$ and $\kappa^{m}(n, \bar{\gamma})=\{\bar{s}\in\mathbb{Z}_{+}^{m}: \,\, \langle\bar{s},\bar{\gamma}\rangle =n\}$, $n\in \mathbb{Z_{+}}$.

\begin{lemma}\label{lem 3}  
{\it Let $\overline{\gamma}^{'} = (\gamma_{1}^{'}, \ldots, \gamma_{m}^{'})$,  $\overline{\gamma} = (\gamma_{1}, \ldots, \gamma_{m})$, $\overline{\theta}=(\theta_{1}, \ldots, \theta_{m})$ and  $0< \gamma_{j}^{'}\leq \gamma_{j}$, $0< \theta_{j}<\infty$, $j=1,\ldots, m$ and $\alpha \in (0, \infty)$. 
Then the following inequality holds:
$$
\left\|\left\{2^{-\alpha\langle\overline{s}, \overline{\gamma}\rangle}
\prod_{j=1}^{m}
(s_{j} + 1)^{\lambda_{j}}
\right\}_{\bar{s}\in Y^{m}(n, \bar{\gamma}^{'})}\right\|_{
l_{\bar \theta}}<<2^{-n\alpha\delta}n^{\sum\limits_{j\in A}\lambda_{j} +\sum\limits_{j\in A\setminus\{j_{1}\}}\frac{1}{\theta_j}},
$$
 where $\delta=\min\{\frac{\gamma_{j}}{\gamma_{j}^{'}} : j=1,\ldots, m\}$, $A=\{j : \frac{\gamma_{j}}{\gamma_{j}^{'}}=\delta, j=1,\ldots, m\}$, $j_{1}=\min\{j : j\in A\}$ and the numbers $\lambda_{j} \in \mathbb{R}$
satisfy the conditions 
$$ 
 \min\{\sum\limits_{j\in A\setminus \{j'\}}\lambda_{j} +\sum\limits_{j\in A\setminus\{j_{1}\}}\frac{1}{\theta_j}, \,\, \lambda_{j'}+\frac{1}{\theta_j'} > 0  
 $$
 and $j' = \max\{j\in A\}$.
}
\end{lemma} 
{\bf Proof.} Let $m=2$ and $1\leq \theta_{j}< \infty$, $j=1, 2$. Then, by the definition of the set
$Y^{2}(n, \bar{\gamma}^{'}) $ we have that 
$$
I_{n}:=\left\|\left\{2^{-\alpha\langle\overline{s}, \overline{\gamma}\rangle}
\prod_{j=1}^{2}(s_{j} + 1)^{\lambda_{j}}
\right\}_{\bar{s}\in Y^{2}(n, \bar{\gamma}^{'})}\right\|_{
l_{\bar \theta}(\mathbb{Z}_{+}^{2})} =
$$
$$
\Bigl\{\sum\limits_{s_{2}< \frac{n}{\gamma_{2}^{'}}}\Bigl[\sum\limits_{s_{1}\geq \frac{n-s_{2}\gamma_{2}^{'}}{\gamma_{1}^{'}}} 2^{-\alpha\langle\overline{s}, \overline{\gamma}\rangle \theta_{1}}\prod_{j=1}^{2}(s_{j} + 1)^{\lambda_{j}\theta_{1}} \Bigr]^{\frac{\theta_{2}}{\theta_{1}}} 
$$
$$
+ \sum\limits_{s_{2}\geq  \frac{n}{\gamma_{2}^{'}}}\Bigl[\sum\limits_{s_{1}=0}^{\infty} 2^{-\alpha\langle\overline{s}, \overline{\gamma}\rangle \theta_{1}}\prod_{j=1}^{2}(s_{j} + 1)^{\lambda_{j}\theta_{1}} \Bigr]^{\frac{\theta_{2}}{\theta_{1}}} \Bigr\}^{\frac{1}{\theta_{2}}}.    \eqno(12)
$$
Since the function $\frac{(t+1)^{\lambda_{j}}}{2^{t\varepsilon}}$ decreases by $ [0, \infty) $, then
$$
\sum\limits_{s_{1}\geq \frac{n-s_{2}\gamma_{2}^{'}}{\gamma_{1}^{'}}} 2^{-s_{1}\gamma_{1}\alpha \theta_{1}}(s_{1} + 1)^{\lambda_{1}\theta_{1}} \leq \Bigl(\frac{n-s_{2}\gamma_{2}^{'}}{\gamma_{1}^{'}}+1\Bigr)^{\lambda_{1}\theta_{1}}2^{-\frac{n-s_{2}\gamma_{2}^{'}}{\gamma_{1}^{'}}\gamma_{1}\varepsilon\theta_{1}}\sum\limits_{s_{1}\geq \frac{n-s_{2}\gamma_{2}^{'}}{\gamma_{1}^{'}}} 2^{-s_{1}\gamma_{1}(\alpha - \varepsilon)\theta_{1}}. \eqno(13)
$$
We choose the number $\varepsilon$ so that $0< \varepsilon < \alpha$. Then from (13), we obtain
$$
\sum\limits_{s_{1}\geq \frac{n-s_{2}\gamma_{2}^{'}}{\gamma_{1}^{'}}} 2^{-s_{1}\gamma_{1}\alpha \theta_{1}}(s_{1} + 1)^{\lambda_{1}\theta_{1}} \leq C 
\Bigl(n-s_{2}\gamma_{2}^{'}+1\Bigr)^{\lambda_{1}\theta_{1}}2^{-\frac{n-s_{2}\gamma_{2}^{'}}{\gamma_{1}^{'}}\gamma_{1}\alpha\theta_{1}}.  \eqno(14)
$$
It can be proved similarly that
$$
\sum\limits_{s_{2}\geq  \frac{n}{\gamma_{2}^{'}}}2^{-s_{2}\gamma_{2}\alpha \theta_{2}}(s_{2} + 1)^{\lambda_{2}\theta_{2}} \leq C(n+1)^{\lambda_{2}\theta_{2}}2^{-n\frac{\gamma_{2}}{\gamma_{2}^{'}}\alpha\theta_{2}}.  \eqno(15)
$$
Now, taking into account that the series 
$$
\sum\limits_{s_{1}=0}^{\infty} 2^{-s_{1}\gamma_{1}\alpha \theta_{1}}(s_{1} + 1)^{\lambda_{1}\theta_{1}}
$$
converges, from inequalities (12), (14), and (15), we obtain.
$$
I_{n} \leq C \left\{2^{-n\frac{\gamma_{1}}{\gamma_{1}^{'}}\alpha \theta_{2}}\sum\limits_{s_{2}< \frac{n}{\gamma_{2}^{'}}} 2^{-s_{2}\gamma_{2}^{'}(\frac{\gamma_{2}}{\gamma_{2}^{'}} - \frac{\gamma_{1}}{\gamma_{1}^{'}})\alpha \theta_{2}}(s_{2}+1)^{\lambda_{2}\theta_{2}} (n-s_{2}\gamma_{2}^{'} + 1)^{\lambda_{1}\theta_{2}}  \right.
$$
$$
\left. + 2^{-n\frac{\gamma_{2}}{\gamma_{2}^{'}}\alpha \theta_{2}}(n + 1)^{\lambda_{2}\theta_{2}}
\right\}^{\frac{1}{\theta_{2}}}.  \eqno(16)
$$
We will estimate the sum
$$
B_{n} = \sum\limits_{s_{2}< \frac{n}{\gamma_{2}^{'}}} 2^{-s_{2}\gamma_{2}^{'}(\frac{\gamma_{2}}{\gamma_{2}^{'}} - \frac{\gamma_{1}}{\gamma_{1}^{'}})\alpha \theta_{2}}(s_{2}+1)^{\lambda_{2}\theta_{2}} (n-s_{2}\gamma_{2}^{'} + 1)^{\lambda_{1}\theta_{2}}.
$$
Let $\frac{\gamma_{2}}{\gamma_{2}^{'}} -  \frac{\gamma_{1}}{\gamma_{1}^{'}} > 0$. Then, by Lemma 2, we have 
$$
B_{n} <<(n+1)^{\lambda_{1}} \eqno (17)
$$
for any  
$\lambda_{1}, \lambda_{2} \in \mathbb{R}$.
It follows from inequalities (16) and (17) that
$$
I_{n} << \left\{2^{-n\frac{\gamma_{1}}{\gamma_{1}^{'}}\alpha \theta_{2}}(n+1)^{\lambda_{1}\theta_{2}} +  2^{-n\frac{\gamma_{2}}{\gamma_{2}^{'}}\alpha \theta_{2}}(n + 1)^{\lambda_{2}\theta_{2}}
\right\}^{\frac{1}{\theta_{2}}}  \eqno(18)
$$
in the case 
$\frac{\gamma_{2}}{\gamma_{2}^{'}} -  \frac{\gamma_{1}}{\gamma_{1}^{'}} > 0$ for  $\lambda_{1}, \lambda_{2} \in \mathbb{R}$. 

Since  
$2^{-n(\frac{\gamma_{2}}{\gamma_{2}^{'}}-\frac{\gamma_{1}}{\gamma_{1}^{'}})\alpha }(n + 1)^{\lambda_{2} - \lambda_{1}} \rightarrow 0$ for $n \rightarrow \infty$, then   
 $$
 2^{-n\frac{\gamma_{2}}{\gamma_{2}^{'}}\alpha}(n + 1)^{\lambda_{2}} \leq 2^{-n\frac{\gamma_{1}}{\gamma_{1}^{'}}\alpha}(n + 1)^{\lambda_{1}} \eqno(19)
 $$
for any natural number $n > n_{0}$, $n_{0} $ is some positive number. 

Now from (18) and (19) it follows that  
$$
I_{n} <<2^{-n\frac{\gamma_{1}}{\gamma_{1}^{'}}\alpha}(n + 1)^{\lambda_{1}} \eqno(20)
$$
in the case 
 $\frac{\gamma_{2}}{\gamma_{2}^{'}} -  \frac{\gamma_{1}}{\gamma_{1}^{'}} > 0$ for  $\lambda_{1}, \lambda_{2} \in \mathbb{R}$. 

Let 
$\frac{\gamma_{2}}{\gamma_{2}^{'}} -  \frac{\gamma_{1}}{\gamma_{1}^{'}} < 0$.
Then, by the second assertion of Lemma 2, we have
$$
B_{n} <<2^{n\bigl(\frac{\gamma_{1}}{\gamma_{1}^{'}}-\frac{\gamma_{2}}{\gamma_{2}^{'}}\bigr)\alpha }(n+1)^{\lambda_{2}} \eqno (21)
$$
for any 
 $\lambda_{1}, \lambda_{2} \in \mathbb{R}$.
It follows from inequalities (16) and (21) that
$$
I_{n} << 2^{-n\frac{\gamma_{2}}{\gamma_{2}^{'}}\alpha}(n + 1)^{\lambda_{2}} \eqno(22)
$$
in the case 
 $\frac{\gamma_{2}}{\gamma_{2}^{'}} -  \frac{\gamma_{1}}{\gamma_{1}^{'}} < 0$ for   $\lambda_{1}, \lambda_{2} \in \mathbb{R}$. 

Let $\frac{\gamma_{2}}{\gamma_{2}^{'}} -  \frac{\gamma_{1}}{\gamma_{1}^{'}} = 0$. 
Then, by Lemma 1, we have    
$$
B_{n} = \sum\limits_{s_{2}< \frac{n}{\gamma_{2}^{'}}} (s_{2}+1)^{\lambda_{2}\theta_{2}} (n-s_{2}\gamma_{2}^{'} + 1)^{\lambda_{1}\theta_{2}} <<(n+1)^{(\lambda_{1}+\lambda_{2})\theta_{2}+1}  \eqno  (23)
$$
if $\min\{\lambda_{1}+ \frac{1}{\theta_{2}}, \,\,  \lambda_{2}+ \frac{1}{\theta_{2}} \} > 0$.

If $\lambda_{2}\theta_{2}+1=0$ and $\lambda_{1}\geq 0$, then 
$$
B_{n} = \sum\limits_{s_{2}< \frac{n}{\gamma_{2}^{'}}} (s_{2}+1)^{-1} (n-s_{2}\gamma_{2}^{'} + 1)^{\lambda_{1}\theta_{2}} <<(n+1)^{\lambda_{1}\theta_{2}}\ln (n+1).   \eqno  (24)
$$

If $\lambda_{2}\theta_{2}+1<0$ and $\lambda_{1}\geq 0$, then 
$$
B_{n}<<(n+1)^{\lambda_{1}\theta_{2}}. 
$$
So if  
  $\lambda_{1}\geq 0$, then
$$
B_{n} << \left\{ \begin{array}{rl} (n+1)^{\lambda_{1}\theta_{2}}, & \mbox{if} \, \,  \lambda_{2}\theta_{2}+1<0 \\
 (n+1)^{\lambda_{1}\theta_{2}}\ln (n+1), & \mbox{if} \, \, \lambda_{2}\theta_{2}+1=0 .
\end{array}
 \right.  \eqno  (25)
$$
Similarly, one can verify that if $\lambda_{2} \geq 0 $, then
$$
B_{n} << \left\{ \begin{array}{rl} (n+1)^{\lambda_{2}\theta_{2}}, & \mbox{if} \, \,  \lambda_{1}\theta_{2}+1<0 \\
 (n+1)^{\lambda_{2}\theta_{2}}\ln (n+1), & \mbox{if} \, \, \lambda_{1}\theta_{2}+1=0 
\end{array}
 \right.  
$$
in the case 
$\frac{\gamma_{2}}{\gamma_{2}^{'}} -  \frac{\gamma_{1}}{\gamma_{1}^{'}} = 0$.

Let $\lambda_{1}< 0$ and $\lambda_{2}< 0$. Then, as in the proof of Lemma 1, we have   
$$
B_{n} = \sum\limits_{s_{2}< \frac{n}{2\gamma_{2}^{'}}} (s_{2}+1)^{\lambda_{2}\theta_{2}} (n-s_{2}\gamma_{2}^{'} + 1)^{\lambda_{1}\theta_{2}} + \sum\limits_{\frac{n}{2\gamma_{2}^{'}}\leq s_{2}< \frac{n}{\gamma_{2}^{'}}} (s_{2}+1)^{\lambda_{2}\theta_{2}} (n-s_{2}\gamma_{2}^{'} + 1)^{\lambda_{1}\theta_{2}}.  \eqno  (26)
$$
It is known that 
$$
\sum\limits_{s_{2}< \frac{n}{2\gamma_{2}^{'}}} (s_{2}+1)^{\lambda_{2}\theta_{2}} << \left\{ \begin{array}{rl} (n+1)^{\lambda_{2}\theta_{2}+1}, & \mbox{if} \, \,  \lambda_{1}\theta_{2}+1>0, \\
 \ln (n+1), & \mbox{if} \, \, \lambda_{2}\theta_{2}+1=0, \\
 1,   & \mbox{if} \, \, \lambda_{2}\theta_{2}+1<0.
\end{array}
 \right.  \eqno  (27)
$$
Similarly 
$$
\sum\limits_{\frac{n}{2\gamma_{2}^{'}}\leq s_{2}< \frac{n}{\gamma_{2}^{'}}} (s_{2}+1)^{\lambda_{2}\theta_{2}} (n-s_{2}\gamma_{2}^{'} + 1)^{\lambda_{1}\theta_{2}} 
<< \left\{ \begin{array}{rl} (n+1)^{\lambda_{1}\theta_{2}+1}, & \mbox{if} \, \,  \lambda_{1}\theta_{2}+1>0, \\
 \ln (n+1), & \mbox{if} \, \, \lambda_{1}\theta_{2}+1=0, \\
 1,   & \mbox{if} \, \, \lambda_{1}\theta_{2}+1<0.
\end{array}
 \right.  \eqno  (28)
$$
Now from inequalities (26)--((29) it follows that
$$
B_{n} = \sum\limits_{s_{2}< \frac{n}{\gamma_{2}^{'}}} (s_{2}+1)^{\lambda_{2}\theta_{2}} (n-s_{2}\gamma_{2}^{'} + 1)^{\lambda_{1}\theta_{2}}
$$
$$
<< \left\{ \begin{array}{rl} (n+1)^{(\lambda_{1}+\lambda_{2})\theta_{2}+1}, & \mbox{if} \, \,  \lambda_{1}\theta_{2}+1>0, \lambda_{2}\theta_{2}+1>0 \\
 (n+1)^{-1}\ln (n+1), & \mbox{if} \, \, \lambda_{1}\theta_{2}+1=\lambda_{2}\theta_{2}+1=0, \\
 \max\{(n+1)^{\lambda_{1}\theta_{2}}, (n+1)^{\lambda_{2}\theta_{2}}\},   & \mbox{if} \, \, \lambda_{1}\theta_{2}+1<0, \lambda_{2}\theta_{2}+1 < 0.
\end{array}
 \right.  \eqno  (29)
$$
in the case 
 $\frac{\gamma_{2}}{\gamma_{2}^{'}} -  \frac{\gamma_{1}}{\gamma_{1}^{'}} = 0$.

From estimates (16), (23), and (29), we obtain
$$
I_{n} << \Bigl\{2^{-n\frac{\gamma_{1}}{\gamma_{1}^{'}}\alpha \theta_{2}}(n+1)^{(\lambda_{1}+\lambda_{2})\theta_{2}+1} +  2^{-n\frac{\gamma_{2}}{\gamma_{2}^{'}}\alpha \theta_{2}}(n + 1)^{\lambda_{2}\theta_{2}}\Bigr\}^{\frac{1}{\theta_{2}}} 
$$
$$
<<2^{-n\alpha\delta}(n+1)^{\lambda_{1}+\lambda_{2}+\frac{1}{\theta_{2}}} \eqno  (30)
$$
in the case 
 $\frac{\gamma_{2}}{\gamma_{2}^{'}} -  \frac{\gamma_{1}}{\gamma_{1}^{'}} = 0$,  if $\min\{\lambda_{1}+\frac{1}{\theta_{2}}, \lambda_{2}+\frac{1}{\theta_{2}}\}> 0$, where $\delta=\frac{\gamma_{1}}{\gamma_{1}^{'}}=\frac{\gamma_{2}}{\gamma_{2}^{'}}$.
If $\lambda_{1}+\frac{1}{\theta_{2}}=0$ and $\lambda_{2}+\frac{1}{\theta_{2}}=0$, then from ( 16) and (29) we obtain
$$
I_{n} <<2^{-n\alpha\delta}\Bigl(\frac{\ln (n+1)}{n+1}\Bigr)^{\frac{1}{\theta_{2}}} \eqno  (31)
$$
in the case  
 $\frac{\gamma_{2}}{\gamma_{2}^{'}} -  \frac{\gamma_{1}}{\gamma_{1}^{'}} = 0$.

If $\lambda_{1}+\frac{1}{\theta_{2}}<0$ and $\lambda_{2}+\frac{1}{\theta_{2}}<0$, then from ( 16) and (29) it follows that
$$
I_{n} <<2^{-n\alpha\delta}\max\{(n+1)^{\lambda_{1}}, (n+1)^{\lambda_{2}}\}
\eqno  (32)
$$
in the case  
$\frac{\gamma_{2}}{\gamma_{2}^{'}} -  \frac{\gamma_{1}}{\gamma_{1}^{'}} = 0$.
This proved the lemma for $1\leq \theta_{j}< \infty$, $j=1,2$.

Let  $\theta_{j} = \infty$, $j=1,2$ and $\overline{s}=(s_{1}, s_{2})\in Y^{2}(n, \overline{\gamma}^{'})$, then  
$$
2^{-\langle \overline{s}, \overline{\gamma} \rangle\alpha}\prod_{j=1}^{2}(s_{j}+1)^{\lambda_{j}} \leq 2^{-n\frac{\gamma_{1}}{\gamma_{1}^{'}}\alpha}2^{-s_{2}\gamma_{2}^{'}(\frac{\gamma_{2}}{\gamma_{2}^{'}}-\frac{\gamma_{1}}{\gamma_{1}^{'}})\alpha}(s_{2}+1)^{\lambda_{2}}(n-s_{2}\gamma_{2}^{'}+1)^{\lambda_{1}}
$$  
for 
 $s_{1}\geq \frac{n-s_{2}\gamma_{2}^{'}}{\gamma_{1}}$ and $0\leq s_{2}<\frac{n}{\gamma_{2}^{'}}$. 

If 
 $s_{1}\geq 0$ and $ s_{2}\geq \frac{n}{\gamma_{2}^{'}}$, then  
$$
2^{-\langle \overline{s}, \overline{\gamma} \rangle\alpha}\prod_{j=1}^{2}(s_{j}+1)^{\lambda_{j}} \leq 2^{-n\frac{\gamma_{2}}{\gamma_{2}^{'}}\alpha}(\frac{n}{\gamma_{2}^{'}}+1)^{\lambda_{2}}.
$$
Let $\frac{\gamma_{2}}{\gamma_{2}^{'}}-\frac{\gamma_{1}}{\gamma_{1}^{'}} > 0$, then choosing the number $\eta\in (0, (\frac{\gamma_{2}}{\gamma_{2}^{'}}-\frac{\gamma_{1}}{\gamma_{1}^{'}})\alpha)$ and considering that the function $2^{t}(t+1)^{\lambda}$ increases by $ [0, + \infty $, $\lambda \in \mathbb{R}$ and $ 0 <n-s_{2}\gamma_{2}^{'} \leq n $, we have.
$$
2^{-\langle \overline{s}, \overline{\gamma} \rangle\alpha}\prod_{j=1}^{2}(s_{j}+1)^{\lambda_{j}}\leq 2^{-n\frac{\gamma_{1}}{\gamma_{1}^{'}}\alpha}2^{-s_{2}\gamma_{2}^{'}(\frac{\gamma_{2}}{\gamma_{2}^{'}}-\frac{\gamma_{1}}{\gamma_{1}^{'}})\alpha}(s_{2}+1)^{\lambda_{2}}(n-s_{2}\gamma_{2}^{'}+1)^{\lambda_{1}} 
$$
$$
 =2^{-n\frac{\gamma_{1}}{\gamma_{1}^{'}}\alpha} 2^{-n\eta}2^{-(n-s_{2}\gamma_{2}^{'})\eta}(n-s_{2}\gamma_{2}^{'}+1)^{\lambda_{1}}2^{-s_{2}\gamma_{2}^{'}((\frac{\gamma_{2}}{\gamma_{2}^{'}}-\frac{\gamma_{1}}{\gamma_{1}^{'}})\alpha-\eta)}(s_{2}+1)^{\lambda_{2}}
$$
$$
\leq 2^{-n\frac{\gamma_{1}}{\gamma_{1}^{'}}\alpha}(n+1)^{\lambda_{1}}
$$
for 
 $s_{1}\geq \frac{n-s_{2}\gamma_{2}^{'}}{\gamma_{1}}$ and $0\leq s_{2}<\frac{n}{\gamma_{2}^{'}}$, $\lambda_{j} \in \mathbb{R}$, $j=1, 2$. 

Let $\frac{\gamma_{2}}{\gamma_{2}^{'}}-\frac{\gamma_{1}}{\gamma_{1}^{'}} < 0$, then choosing the number $\eta\in (0, -(\frac{\gamma_{2}}{\gamma_{2}^{'}}-\frac{\gamma_{1}}{\gamma_{1}^{'}})\alpha)$, similarly one can prove that
$$
2^{-\langle \overline{s}, \overline{\gamma} \rangle\alpha}\prod_{j=1}^{2}(s_{j}+1)^{\lambda_{j}}\leq 2^{-n\frac{\gamma_{2}}{\gamma_{2}^{'}}\alpha}(n+1)^{\lambda_{2}}
$$
for 
 $s_{1}\geq \frac{n-s_{2}\gamma_{2}^{'}}{\gamma_{1}}$ and $0\leq s_{2}<\frac{n}{\gamma_{2}^{'}}$, $\lambda_{j} \in \mathbb{R}$, $j=1, 2$. 

If 
 $\frac{\gamma_{2}}{\gamma_{2}^{'}}-\frac{\gamma_{1}}{\gamma_{1}^{'}}=0$, then 
$$
2^{-\langle \overline{s}, \overline{\gamma} \rangle\alpha}\prod_{j=1}^{2}(s_{j}+1)^{\lambda_{j}}\leq 2^{-n\frac{\gamma_{1}}{\gamma_{1}^{'}}\alpha}(n+1)^{\lambda_{2}}(s_{2}+1)^{\lambda_{2}}(n-s_{2}\gamma_{2}^{'}+1)^{\lambda_{1}} \leq 2^{-n\frac{\gamma_{1}}{\gamma_{1}^{'}}\alpha}(n+1)^{\lambda_{1} + \lambda_{2}}
$$
for 
 $s_{1}\geq \frac{n-s_{2}\gamma_{2}^{'}}{\gamma_{1}}$ and $0\leq s_{2}<\frac{n}{\gamma_{2}^{'}}$, $\lambda_{j} \geq 0$, $j=1, 2$. 

Now, it follows from these inequalities that
$$
I_{n}=\sup_{\overline{s}\in Y^{2}(n, \overline{\gamma}^{'})}2^{-\langle \overline{s}, \overline{\gamma} \rangle\alpha}\prod_{j=1}^{2}(s_{j}+1)^{\lambda_{j}} \leq 2^{-n\alpha\delta}(n+1)^{\sum\limits_{j\in A} \lambda_{j}}, 
$$
where 
 $\lambda_{j}\in \mathbb{R}$, $j=1, 2$ if 
  $A\setminus \{j_{1}\}=\emptyset$ and $\lambda_{j}\geq 0$, $j=1, 2$ if  $A\setminus \{j_{1}\}\neq \emptyset$, $j_{1}=\min\{j\in A\}$. 
This proves the lemma for $ m = 2 $.

Now consider the case $m=3$.
$$
I_{n}:=\left\|\left\{2^{-\alpha\langle\overline{s}, \overline{\gamma}\rangle}
\prod_{j=1}^{3}(s_{j} + 1)^{\lambda_{j}}
\right\}_{\bar{s}\in Y^{3}(n, \bar{\gamma}^{'})}\right\|_{
l_{\bar \theta}(\mathbb{Z}_{+}^{3})} \leq
$$
$$
\leq \Biggl[\sum\limits_{s_{3}< \frac{n}{\gamma_{3}^{'}}}\left(2^{-s_{3}\gamma_{3}\alpha} (s_{3}+1)^{\lambda_{3}}\left\|\left\{2^{-\alpha\langle\overline{s}_{2}, \overline{\gamma}_{2}\rangle}
\prod_{j=1}^{3}(s_{j} + 1)^{\lambda_{j}}
\right\}_{\bar{s}_{2}\in Y^{2}(n - s_{3}\gamma_{3}^{'}, \bar{\gamma}^{'})}\right\|_{
l_{\bar\theta_{2}}(\mathbb{Z}_{+}^{2})} \right)^{\theta_{3}}\Biggr]^{\frac{1}{\theta_{3}}} + 
$$
$$
\Biggl[\sum\limits_{s_{3}\geq \frac{n}{\gamma_{3}^{'}}}  \left(2^{-s_{3}\gamma_{3}\alpha} (s_{3}+1)^{\lambda_{3}}\left\|\left\{2^{-\alpha\langle\overline{s}_{2}, \overline{\gamma}_{2}\rangle}
\prod_{j=1}^{3}(s_{j} + 1)^{\lambda_{j}}
\right\}_{\bar{s}_{2}\in \mathbb{Z}_{+}^{2}}\right\|_{
l_{\bar\theta_{2}(\mathbb{Z}_{+}^{2})}} \right)^{\theta_{3}}\Biggr]^{\frac{1}{\theta_{3}}} 
$$
$$
=\sigma_{1}(n) + \sigma_{2}(n),   \eqno  (33) 
$$
where  $\bar{a}_{2}=(a_{1}, a_{2})$.

If $\frac{\gamma_{2}}{\gamma_{2}^{'}}-\frac{\gamma_{1}}{\gamma_{1}^{'}}> 0$ and $ \lambda_{1}, \lambda_{2}\in \mathbb{R}$, then according to estimate (20) with $ n $ replaced by $ n - s_{3} \gamma_{3}^{'}$ we have
$$
I_{n- s_{3}\gamma_{3}^{'}} = \left\|\left\{2^{-\alpha\langle\overline{s}_{2}, \overline{\gamma}_{2}\rangle}
\prod_{j=1}^{3}(s_{j} + 1)^{\lambda_{j}}
\right\}_{\bar{s}_{2}\in Y^{2}(n - s_{3}\gamma_{3}^{'}, \bar{\gamma}^{'})}\right\|_{
l_{\bar\theta_{2}(\mathbb{Z}_{+}^{2})}} 
$$
$$
<< 2^{-(n-s_{3}\gamma_{3}^{'})\frac{\gamma_{1}}{\gamma_{1}^{'}}\alpha}(n - s_{3}\gamma_{3}^{'} + 1)^{\lambda_{1}} \eqno(34)
$$
for  $0\leq s_{3} <\frac{n}{\gamma_{3}^{'}}$.

According to estimate (34), we have
$$
\sigma_{1}(n) <<
2^{-n\frac{\gamma_{1}}{\gamma_{1}^{'}}\alpha}\Bigl[\sum\limits_{s_{3}< \frac{n}{\gamma_{3}^{'}}} 2^{-s_{3}\gamma_{3}^{'}(\frac{\gamma_{3}}{\gamma_{3}^{'}} - \frac{\gamma_{1}}{\gamma_{1}^{'}})\alpha \theta_{3}}(s_{3}+1)^{\lambda_{3}\theta_{3}} (n-s_{3}\gamma_{3}^{'} + 1)^{\lambda_{1}\theta_{3}}\Bigr]^{\frac{1}{\theta_{3}}}.  \eqno(35)
$$

Let $\frac{\gamma_{3}}{\gamma_{3}^{'}} - \frac{\gamma_{1}}{\gamma_{1}^{'}} > 0$.
If $\frac{\gamma_{2}}{\gamma_{2}^{'}} - \frac{\gamma_{1}}{\gamma_{1}^{'}} > 0$, then by the lemma 2 we get
$$
\sum\limits_{s_{3}< \frac{n}{\gamma_{3}^{'}}} 2^{-s_{3}\gamma_{3}^{'}(\frac{\gamma_{3}}{\gamma_{3}^{'}} - \frac{\gamma_{1}}{\gamma_{1}^{'}})\alpha \theta_{3}}(s_{3}+1)^{\lambda_{3}\theta_{3}} (n-s_{3}\gamma_{3}^{'} + 1)^{\lambda_{1}\theta_{3}} <<(n+1)^{\lambda_{1}\theta_{3}}   \eqno(36)
$$
for 
$\lambda_{j}\in \mathbb{R}$, $j=1, 2, 3$.
It now follows from inequalities (35) and (36) that
$$
\sigma_{1}(n)<<
2^{-n\frac{\gamma_{1}}{\gamma_{1}^{'}}\alpha}(n+1)^{\lambda_{1}} \eqno(37)
$$
in the case  
 $\frac{\gamma_{3}}{\gamma_{3}^{'}} - \frac{\gamma_{1}}{\gamma_{1}^{'}} > 0$ and $\frac{\gamma_{2}}{\gamma_{2}^{'}} - \frac{\gamma_{1}}{\gamma_{1}^{'}} > 0$ ,  for  $\delta = \frac{\gamma_{1}}{\gamma_{1}^{'}}$.

If $\frac{\gamma_{2}}{\gamma_{2}^{'}} - \frac{\gamma_{1}}{\gamma_{1}^{'}} < 0$, then according to what was proved inequality (22) with $n$ replaced by $n - s_{3} \gamma_{3}^{'}$ we have
$$
I_{n- s_{3}\gamma_{3}^{'}} = \left\|\left\{2^{-\alpha\langle\overline{s}_{2}, \overline{\gamma}_{2}\rangle}
\prod_{j=1}^{3}(s_{j} + 1)^{\lambda_{j}}
\right\}_{\bar{s}_{2}\in Y^{2}(n - s_{3}\gamma_{3}^{'}, \bar{\gamma}^{'})}\right\|_{
l_{\bar\theta_{2}(\mathbb{Z}_{+}^{2})}} 
$$
$$
<<2^{-(n-s_{3}\gamma_{3}^{'})\frac{\gamma_{2}}{\gamma_{2}^{'}}\alpha}(n - s_{3}\gamma_{3}^{'} + 1)^{\lambda_{2}} \eqno(38)
$$
for 
 $0\leq s_{3} <\frac{n}{\gamma_{3}^{'}}$.
Therefore (see (35))
$$
\sigma_{1}(n) <<
2^{-n\frac{\gamma_{2}}{\gamma_{2}^{'}}\alpha}\Bigl[\sum\limits_{s_{3}< \frac{n}{\gamma_{3}^{'}}} 2^{-s_{3}\gamma_{3}^{'}(\frac{\gamma_{3}}{\gamma_{3}^{'}} - \frac{\gamma_{2}}{\gamma_{2}^{'}})\alpha \theta_{3}}(s_{3}+1)^{\lambda_{3}\theta_{3}} (n-s_{3}\gamma_{3}^{'} + 1)^{\lambda_{2}\theta_{3}}\Bigr]^{\frac{1}{\theta_{3}}}.  \eqno(39)
$$
If 
 $\frac{\gamma_{3}}{\gamma_{3}^{'}} - \frac{\gamma_{1}}{\gamma_{1}^{'}} > 0$ and $\frac{\gamma_{2}}{\gamma_{2}^{'}} - \frac{\gamma_{1}}{\gamma_{1}^{'}} < 0$, then $\frac{\gamma_{3}}{\gamma_{3}^{'}} - \frac{\gamma_{2}}{\gamma_{2}^{'}} > 0$. 
Hence, by Lemma 2, from (39) we obtain 
 $$
\sigma_{1}(n) <<2^{-n\frac{\gamma_{2}}{\gamma_{2}^{'}}\alpha}(n+1)^{\lambda_{2}}  \eqno(40)
$$
in the case 
 $\frac{\gamma_{3}}{\gamma_{3}^{'}} - \frac{\gamma_{1}}{\gamma_{1}^{'}} > 0$ and $\frac{\gamma_{2}}{\gamma_{2}^{'}} - \frac{\gamma_{1}}{\gamma_{1}^{'}} < 0$, for $\delta = \frac{\gamma_{2}}{\gamma_{2}^{'}}$ and for $\lambda_{j}\in \mathbb{R}$, $j=1, 2, 3$.

Estimates (37) and (40) can be jointly written in the following form
$$
\sigma_{1}(n) <<2^{-n\alpha\delta}(n+1)^{\sum\limits_{j\in A}\lambda_{j} + \sum\limits_{j\in A\setminus \{j_{0}\}} \frac{1}{\theta_{j}}}  \eqno(41)
$$
in the case 
 $\frac{\gamma_{3}}{\gamma_{3}^{'}} - \frac{\gamma_{1}}{\gamma_{1}^{'}} > 0$, where $\delta = \min\{\frac{\gamma_{j}}{\gamma_{j}^{'}} : \,\, j=1,2,3\}$ and $j_{0}=\min\{j=1, 2, 3: \,\, \delta =\frac{\gamma_{j}}{\gamma_{j}^{'}}\}$.

Let  
 $\frac{\gamma_{3}}{\gamma_{3}^{'}} - \frac{\gamma_{1}}{\gamma_{1}^{'}} < 0$.
 If $\frac{\gamma_{2}}{\gamma_{2}^{'}} - \frac{\gamma_{1}}{\gamma_{1}^{'}} > 0$, then from the formula (35) by Lemma 2 we obtain
$$
\sigma_{1}(n) <<2^{-n\frac{\gamma_{3}}{\gamma_{3}^{'}}\alpha}(n+1)^{\lambda_{3}}  \eqno(42)
$$
for 
 $\lambda_{j}\in \mathbb{R}$, $j=1, 2, 3$.

Let 
 $\frac{\gamma_{3}}{\gamma_{3}^{'}} - \frac{\gamma_{1}}{\gamma_{1}^{'}} < 0$ and $\frac{\gamma_{2}}{\gamma_{2}^{'}} - \frac{\gamma_{1}}{\gamma_{1}^{'}} < 0$. 
Since $\frac{\gamma_{2}}{\gamma_{2}^{'}} - \frac{\gamma_{1}}{\gamma_{1}^{'}} < 0$, then in formula (22) replacing $n$ by $n - s_{3}\gamma_{3}^{'}$ we obtain 
$$
I_{n- s_{3}\gamma_{3}^{'}} <<2^{-(n-s_{3}\gamma_{3}^{'})\frac{\gamma_{2}}{\gamma_{2}^{'}}\alpha}(n - s_{3}\gamma_{3}^{'} + 1)^{\lambda_{2}}  \eqno(43)
$$
for 
$\lambda_{j}\in \mathbb{R}$, $j=1, 2$.  Therefore 
$$
\sigma_{1}(n)<<
2^{-n\frac{\gamma_{2}}{\gamma_{2}^{'}}\alpha}\Bigl[\sum\limits_{s_{3}< \frac{n}{\gamma_{3}^{'}}} 2^{s_{3}\gamma_{3}^{'}(\frac{\gamma_{2}}{\gamma_{2}^{'}} - \frac{\gamma_{3}}{\gamma_{3}^{'}})\alpha \theta_{3}}(s_{3}+1)^{\lambda_{3}\theta_{3}} (n-s_{3}\gamma_{3}^{'} + 1)^{\lambda_{2}\theta_{3}}\Bigr]^{\frac{1}{\theta_{3}}}  \eqno(44)
$$
in the case  
$\frac{\gamma_{3}}{\gamma_{3}^{'}} - \frac{\gamma_{1}}{\gamma_{1}^{'}} < 0$ and $\frac{\gamma_{2}}{\gamma_{2}^{'}} - \frac{\gamma_{1}}{\gamma_{1}^{'}} < 0$, 
for  
 $\lambda_{j}\in \mathbb{R}$, $j=1, 2, 3$.
If 
$\frac{\gamma_{2}}{\gamma_{2}^{'}} - \frac{\gamma_{3}}{\gamma_{3}^{'}} > 0$, then by Lemma 2 from (44) we obtain
$$
\sigma_{1}(n) <<2^{-n\frac{\gamma_{3}}{\gamma_{3}^{'}}\alpha}(n+1)^{\lambda_{3}}  \eqno(45)
$$
in the case 
 $\frac{\gamma_{3}}{\gamma_{3}^{'}} - \frac{\gamma_{1}}{\gamma_{1}^{'}} < 0$ and $\frac{\gamma_{2}}{\gamma_{2}^{'}} - \frac{\gamma_{1}}{\gamma_{1}^{'}} < 0$
for 
 $\lambda_{j}\in \mathbb{R}$, $j=1, 2, 3$.

Let $\frac{\gamma_{3}}{\gamma_{3}^{'}} - \frac{\gamma_{1}}{\gamma_{1}^{'}} < 0$ and $\frac{\gamma_{2}}{\gamma_{2}^{'}} - \frac{\gamma_{1}}{\gamma_{1}^{'}} < 0$, but
$\frac{\gamma_{2}}{\gamma_{2}^{'}} - \frac{\gamma_{3}}{\gamma_{3}^{'}} < 0$. Then, according to Lemma 2, from (44) we obtain
$$
\sigma_{1}(n) <<2^{-n\frac{\gamma_{2}}{\gamma_{2}^{'}}\alpha}(n+1)^{\lambda_{2}}  \eqno(46)
$$
in the case  
 $\frac{\gamma_{3}}{\gamma_{3}^{'}} - \frac{\gamma_{1}}{\gamma_{1}^{'}} < 0$ and $\frac{\gamma_{2}}{\gamma_{2}^{'}} - \frac{\gamma_{1}}{\gamma_{1}^{'}} < 0$
for  $\lambda_{j}\in \mathbb{R}$, $j=1, 2, 3$.

Now estimates (42), (43), (45), and (46) can be jointly written in the following form
$$
\sigma_{1}(n)<<2^{-n\alpha\delta}(n+1)^{\sum\limits_{j\in A}\lambda_{j} + \sum\limits_{j\in A\setminus \{j_{0}\}} \frac{1}{\theta_{j}}}  \eqno(47)
$$
in the case  
 $\frac{\gamma_{3}}{\gamma_{3}^{'}} - \frac{\gamma_{1}}{\gamma_{1}^{'}} < 0$, where 
  $\delta = \min\{\frac{\gamma_{j}}{\gamma_{j}^{'}} : \,\, j=1,2,3\}$ and $j_{0}=\min\{j=1, 2, 3: \,\, \delta =\frac{\gamma_{j}}{\gamma_{j}^{'}}\}$.

Let  
 $\frac{\gamma_{3}}{\gamma_{3}^{'}} - \frac{\gamma_{1}}{\gamma_{1}^{'}} = 0$.
Then  (see (34) and (35))  
$$
\sigma_{1}(n) <<
2^{-n\frac{\gamma_{1}}{\gamma_{1}^{'}}\alpha}\Bigl[\sum\limits_{s_{3}< \frac{n}{\gamma_{3}^{'}}} (s_{3}+1)^{\lambda_{3}\theta_{3}} (n-s_{3}\gamma_{3}^{'} + 1)^{\lambda_{1}\theta_{3}}\Bigr]^{\frac{1}{\theta_{3}}} 
$$
if 
 $\frac{\gamma_{2}}{\gamma_{2}^{'}} - \frac{\gamma_{1}}{\gamma_{1}^{'}} > 0$,
for
 $\lambda_{j}\in \mathbb{R}$, $j=1, 2, 3$. According to Lemma 1, from this we obtain    
$$
\sigma_{1}(n) <<
2^{-n\frac{\gamma_{1}}{\gamma_{1}^{'}}\alpha}(n+1)^{\lambda_{1}+\lambda_{3}+\frac{1}{\theta_{3}}} \eqno(48)
$$
for $\min\{\lambda_{1}+\frac{1}{\theta_{3}}, \lambda_{3}+\frac{1}{\theta_{3}}\} > 0$ and $\lambda_{2}\in \mathbb{R}$.

Let 
 $\frac{\gamma_{3}}{\gamma_{3}^{'}} - \frac{\gamma_{1}}{\gamma_{1}^{'}} = 0$ and $\frac{\gamma_{2}}{\gamma_{2}^{'}} - \frac{\gamma_{1}}{\gamma_{1}^{'}} < 0$.
Then (see (39)) 
$$
\sigma_{1}(n) <<2^{-n\frac{\gamma_{2}}{\gamma_{2}^{'}}\alpha}\Bigl[\sum\limits_{s_{3}< \frac{n}{\gamma_{3}^{'}}} 2^{-s_{3}\gamma_{3}^{'}(\frac{\gamma_{3}}{\gamma_{3}^{'}} - \frac{\gamma_{2}}{\gamma_{2}^{'}})\alpha \theta_{3}}(s_{3}+1)^{\lambda_{3}\theta_{3}} (n-s_{3}\gamma_{3}^{'} + 1)^{\lambda_{2}\theta_{3}}\Bigr]^{\frac{1}{\theta_{3}}}.  \eqno(49)
$$
Taking into account that $\frac{\gamma_{3}}{\gamma_{3}^{'}} - \frac{\gamma_{1}}{\gamma_{1}^{'}} = 0$ we have $\frac{\gamma_{3}}{\gamma_{3}^{'}} - \frac{\gamma_{2}}{\gamma_{2}^{'}} = \frac{\gamma_{1}}{\gamma_{1}^{'}} - \frac{\gamma_{2}}{\gamma_{2}^{'}} > 0$. Therefore, according to Lemma 2, from (49) we obtain
$$
\sigma_{1}(n)<<
2^{-n\frac{\gamma_{2}}{\gamma_{2}^{'}}\alpha}(n+1)^{\lambda_{2}}= C
2^{-n\alpha\delta}(n+1)^{\lambda_{2}} \eqno(50)
$$
for $\lambda_{j}\in \mathbb{R}$, $j=1, 2, 3$,   in the case 
 $\frac{\gamma_{3}}{\gamma_{3}^{'}} - \frac{\gamma_{1}}{\gamma_{1}^{'}} = 0$ and $\frac{\gamma_{2}}{\gamma_{2}^{'}} - \frac{\gamma_{1}}{\gamma_{1}^{'}} < 0$,  $\delta=\frac{\gamma_{2}}{\gamma_{2}^{'}}$ and $A=\{2\}$.

Let 
 $\frac{\gamma_{3}}{\gamma_{3}^{'}} - \frac{\gamma_{1}}{\gamma_{1}^{'}} = 0$ and $\frac{\gamma_{2}}{\gamma_{2}^{'}} - \frac{\gamma_{1}}{\gamma_{1}^{'}} = 0$. Since   
$\frac{\gamma_{2}}{\gamma_{2}^{'}} - \frac{\gamma_{1}}{\gamma_{1}^{'}} = 0$, then 
$$
I_{n- s_{3}\gamma_{3}^{'}}<<2^{-(n-s_{3}\gamma_{3}^{'})\alpha\delta}(n - s_{3}\gamma_{3}^{'} + 1)^{\lambda_{1}+\lambda_{2}+\frac{1}{\theta_{2}}}  
$$
for $\min\{\lambda_{1}+\frac{1}{\theta_{2}}, \lambda_{2}+\frac{1}{\theta_{2}}\} > 0$. Therefore, since $\frac{\gamma_{3}}{\gamma_{3}^{'}} - \frac{\gamma_{1}}{\gamma_{1}^{'}} = 0$, then
$$
\sigma_{1}(n)<<
2^{-n\frac{\gamma_{1}}{\gamma_{1}^{'}}\alpha}\Bigl[\sum\limits_{s_{3}< \frac{n}{\gamma_{3}^{'}}} (s_{3}+1)^{\lambda_{3}\theta_{3}} (n-s_{3}\gamma_{3}^{'} + 1)^{(\lambda_{1}+ \lambda_{2}+\frac{1}{\theta_{2}})\theta_{3}}\Bigr]^{\frac{1}{\theta_{3}}}. 
$$
Now, according to Lemma 1, from this we obtain
$$
\sigma_{1}(n)<<
2^{-n\frac{\gamma_{1}}{\gamma_{1}^{'}}\alpha}(n-s_{3}\gamma_{3}^{'} + 1)^{\lambda_{1}+ \lambda_{2}+ \lambda_{3}+\frac{1}{\theta_{2}} + \frac{1}{\theta_{3}}}   \eqno(51)
 $$
for $\min\{\lambda_{1}+\lambda_{2}+\frac{1}{\theta_{2}}+\frac{1}{\theta_{3}}, \lambda_{3}+\frac{1}{\theta_{3}}\} > 0$,  
$\frac{\gamma_{3}}{\gamma_{3}^{'}} - \frac{\gamma_{1}}{\gamma_{1}^{'}} = 0$ and $\frac{\gamma_{2}}{\gamma_{2}^{'}} - \frac{\gamma_{1}}{\gamma_{1}^{'}} = 0$.

If the set $A\setminus \{j_{0}\}$ is not empty, then estimates (48) and (51) can be written together in the following form
$$
\sigma_{1}(n)<<2^{-n\alpha\delta}(n+1)^{\sum\limits_{j\in A}\lambda_{j} + \sum\limits_{j\in A\setminus \{j_{0}\}} \frac{1}{\theta_{j}}}  \eqno(52)
$$
for 
 $\min\{\sum\limits_{j\in A\setminus \{j'\}}\lambda_{j} + \sum\limits_{j\in A\setminus \{j_{0}\}} \frac{1}{\theta_{j}}, \lambda_{3}+\frac{1}{\theta_{3}}\} > 0$, 
where $\delta = \min\{\frac{\gamma_{j}}{\gamma_{j}^{'}} : \,\, j=1,2,3\}$ and $j_{0}=\min\{j=1, 2, 3: \,\, \delta =\frac{\gamma_{j}}{\gamma_{j}^{'}}\}$,  $j^{'}=\max\{j\in A\}$.

Now estimates (41), (47), and (52) can be written together in the following form
 $$
\sigma_{1}(n)<<2^{-n\alpha\delta}(n+1)^{\sum\limits_{j\in A}\lambda_{j} + \sum\limits_{j\in A\setminus \{j_{0}\}} \frac{1}{\theta_{j}}}  \eqno(53)
$$
for 
 $\min\{\sum\limits_{j\in A\setminus \{j'\}}\lambda_{j} + \sum\limits_{j\in A\setminus \{j_{0}\}} \frac{1}{\theta_{j}}, \lambda_{j^{'}}+\frac{1}{\theta_{j^{'}}}\} > 0$, 
$\delta = \min\{\frac{\gamma_{j}}{\gamma_{j}^{'}} : \,\, j=1,2,3\}$ and $j_{0}=\min\{j=1, 2, 3: \,\, \delta =\frac{\gamma_{j}}{\gamma_{j}^{'}}\}$,  $j^{'}=\max\{j\in A\}$, $A=\{j=1, 2, 3: \delta =\frac{\gamma_{j}}{\gamma_{j}^{'}}\}$.

Further, it follows from inequalities (33) and (53) that
$$
I_{n} <<\Bigl[2^{-n\alpha\delta}(n+1)^{\sum\limits_{j\in A}\lambda_{j} + \sum\limits_{j\in A\setminus \{j_{0}\}} \frac{1}{\theta_{j}}}  
$$
$$
+
\Biggl[\sum\limits_{s_{3}\geq \frac{n}{\gamma_{3}^{'}}}  \left(2^{-s_{3}\gamma_{3}\alpha} (s_{3}+1)^{\lambda_{3}}\left\|\left\{2^{-\alpha\langle\overline{s}_{2}, \overline{\gamma}_{2}\rangle}
\prod_{j=1}^{3}(s_{j} + 1)^{\lambda_{j}}
\right\}_{\bar{s}_{2}\in \mathbb{Z}_{+}^{2}}\right\|_{
l_{\bar\theta_{2}(\mathbb{Z}_{+}^{2})}} \right)^{\theta_{3}}\Biggr]^{\frac{1}{\theta_{3}}} \Bigr] \eqno(54)
$$
for 
$\min\{\sum\limits_{j\in A\setminus \{j'\}}\lambda_{j} + \sum\limits_{j\in A\setminus \{j_{0}\}} \frac{1}{\theta_{j}}, \lambda_{j^{'}}+\frac{1}{\theta_{j^{'}}}\} > 0$, 
$\delta = \min\{\frac{\gamma_{j}}{\gamma_{j}^{'}} : \,\, j=1,2,3\}$ and $j_{0}=\min\{j=1, 2, 3: \,\, \delta =\frac{\gamma_{j}}{\gamma_{j}^{'}}\}$,  $j^{'}=\max\{j\in A\}$, $A=\{j=1, 2, 3: \delta =\frac{\gamma_{j}}{\gamma_{j}^{'}}\}$.

Since 
 $\alpha>0$, then 
$$
\left\|\left\{2^{-\alpha\langle\overline{s}_{2}, \overline{\gamma}_{2}\rangle}
\prod_{j=1}^{3}(s_{j} + 1)^{\lambda_{j}}
\right\}_{\bar{s}_{2}\in \mathbb{Z}_{+}^{2}}\right\|_{
l_{\bar\theta_{2}(\mathbb{Z}_{+}^{2})}} < \infty
$$
and 
$$
\sum\limits_{s_{3}\geq \frac{n}{\gamma_{3}^{'}}}  2^{-s_{3}\gamma_{3}\alpha\theta_{3}} (s_{3}+1)^{\lambda_{3}\theta_{3}} \leq C 2^{-n\frac{\gamma_{3}}{\gamma_{3}^{'}}\alpha\theta_{3}}(n+1)^{\lambda_{3}\theta_{3}}.
$$
Therefore, from (54) we obtain
$$
I_{n}<<\Bigl[2^{-n\alpha\delta}(n+1)^{\sum\limits_{j\in A}\lambda_{j} + \sum\limits_{j\in A\setminus \{j_{0}\}} \frac{1}{\theta_{j}}} + 2^{-n\frac{\gamma_{3}}{\gamma_{3}^{'}}\alpha\theta_{3}}(n+1)^{\lambda_{3}\theta_{3}}\Bigr] \eqno(55)
 $$
for  
$\min\{\sum\limits_{j\in A\setminus \{j'\}}\lambda_{j} + \sum\limits_{j\in A\setminus \{j_{0}\}} \frac{1}{\theta_{j}}, \lambda_{j^{'}}+\frac{1}{\theta_{j^{'}}}\} > 0$, 
$\delta = \min\{\frac{\gamma_{j}}{\gamma_{j}^{'}} : \,\, j=1,2,3\}$ and $j_{0}=\min\{j=1, 2, 3: \,\, \delta =\frac{\gamma_{j}}{\gamma_{j}^{'}}\}$,  $j^{'}=\max\{j\in A\}$, $A=\{j=1, 2, 3: \delta =\frac{\gamma_{j}}{\gamma_{j}^{'}}\}$.

Since 
 $\min\{\sum\limits_{j\in A\setminus \{j'\}}\lambda_{j} + \sum\limits_{j\in A\setminus \{j_{0}\}} \frac{1}{\theta_{j}}, \lambda_{j^{'}}+\frac{1}{\theta_{j^{'}}}\} > 0$, then   
$$
2^{-n\frac{\gamma_{3}}{\gamma_{3}^{'}}\alpha\theta_{3}}(n+1)^{\lambda_{3}\theta_{3}}<<2^{-n\alpha\delta}(n+1)^{\sum\limits_{j\in A}\lambda_{j} + \sum\limits_{j\in A\setminus \{j_{0}\}} \frac{1}{\theta_{j}}}.  \eqno(56)
$$
It now follows from inequalities (55) and (56) that  
$$
I_{n} <<2^{-n\alpha\delta}(n+1)^{\sum\limits_{j\in A}\lambda_{j}} \eqno (57)  
$$
for 
 $\min\{\sum\limits_{j\in A\setminus \{j'\}}\lambda_{j} + \sum\limits_{j\in A\setminus \{j_{0}\}} \frac{1}{\theta_{j}}, \lambda_{j^{'}}+\frac{1}{\theta_{j^{'}}}\} > 0$, 
$\delta = \min\{\frac{\gamma_{j}}{\gamma_{j}^{'}} : \,\, j=1,2,3\}$ and $j_{0}=\min\{j=1, 2, 3: \,\, \delta =\frac{\gamma_{j}}{\gamma_{j}^{'}}\}$,  $j^{'}=\max\{j\in A\}$, $A=\{j=1, 2, 3: \delta =\frac{\gamma_{j}}{\gamma_{j}^{'}}\}$.
This proves the lemma for $ m = 3 $.

Suppose that the assertion of the lemma is true for $m-1 \geq 2 $, that is,
$$
\left\|\left\{2^{-\alpha\langle\overline{s}_{m-1}, \overline{\gamma}_{m-1}\rangle}
\prod_{j=1}^{m-1}
(s_{j}+1)^{\lambda_{j}}
\right\}_{\bar{s}_{m-1}\in Y^{m-1}(n, \bar{\gamma}^{'})}\right\|_{
l_{\bar \theta}(\mathbb{Z}^{m-1}_{+})}
$$
$$
<<2^{-n\alpha\delta_{m-1}}(n+1)^{\sum\limits_{j\in A_{m-1}}\lambda_{j}+\sum\limits_{j\in A_{m-1}\setminus\{j_{0}\}}\frac{1}{\theta_j}},
\eqno (58) 
$$
under the condition 
 $\min\{\sum\limits_{j\in A_{m-1}\setminus \{j'\}}\lambda_{j} + \sum\limits_{j\in A_{m-1}\setminus \{j_{0}\}} \frac{1}{\theta_{j}}, \lambda_{j^{'}}+\frac{1}{\theta_{j^{'}}}\} > 0$,
where 
 $\delta_{m-1}=\min\{\frac{\gamma_{j}}{\gamma_{j}^{'}} : j=1,\ldots, m-1\}$, $A_{m-1}=\{j : \frac{\gamma_{j}}{\gamma_{j}^{'}}=\delta, j=1,\ldots, m-1\}$, $j_{0}=\min\{j : j\in A_{m-1}\}$, $j^{'}=\max\{j\in A_{m-1}\}$, $\overline{a}_{m-1}=(a_{1},...,a_{m-1})$.

We prove the lemma for $m$. By the definition of the set $Y^{m}(n, \bar{\gamma}^{'})$ and by the Minkowski inequality, we obtain
$$
I_{n}=\left\|\left\{2^{-\alpha\langle\overline{s}, \overline{\gamma}\rangle}
\prod_{j=1}^{m}(s_{j}+1)^{\lambda_{j}}
\right\}_{\bar{s}\in Y^{m}(n, \bar{\gamma}^{'})}\right\|_{
l_{\bar \theta}}<<
$$
$$
 \left\{\left[\sum\limits_{0\leq s_{m}< \frac{n}{\gamma_{m}^{'}}}\Bigl(2^{-s_{m}\gamma_{m}\alpha}(s_{m}+1)^{\lambda_{m}} \left\|\left\{2^{-\alpha\langle\overline{s}, \overline{\gamma}\rangle}
\prod_{j=1}^{m-1}(s_{j}+1)^{\lambda_{j}}
\right\}_{\bar{s}\in Y^{m-1}(n-s_{m}\gamma_{m}, \bar{\gamma}^{'})}\right\|_{
l_{\bar \theta}(\mathbb{Z}_{+}^{m-1})} \Bigr)^{\theta_{m}} \right]^{\frac{1}{\theta_{m}}} \right. 
$$
$$
\left. +  \left[\sum\limits_{ s_{m}\geq  \frac{n}{\gamma_{m}^{'}}}\Bigl(2^{-s_{m}\gamma_{m}\alpha}(s_{m}+1)^{\lambda_{m}}\left\|\left\{2^{-\alpha\langle\overline{s}, \overline{\gamma}\rangle}
\prod_{j=1}^{m-1}(s_{j}+1)^{\lambda_{j}}
\right\}_{\bar{s}\in \mathbb{Z}_{+}^{m-1}}\right\|_{
l_{\bar \theta}(\mathbb{Z}_{+}^{m-1})} \Bigr)^{\theta_{m}} \right]^{\frac{1}{\theta_{m}}}\right\} 
$$
$$
= C\left\{\sigma_{1}(n) + \sigma_{2}(n)\right\}. \eqno (59) 
$$
In order to estimate $\sigma_{1}(n)$, we use inequality (58). Then
$$
\sigma_{1}(n) <<\left[\sum\limits_{0\leq s_{m}< \frac{n}{\gamma_{m}^{'}}}\Bigl(2^{-s_{m}\gamma_{m}\alpha}(s_{m}+1)^{\lambda_{m}} \right.
$$
$$
\left.
\times 2^{-(n-s_{m}\gamma_{m}^{'})\alpha\delta_{m-1}}\prod_{j=1}^{m-1}(s_{j}+1)^{\lambda_{j}}(n-s_{m}\gamma_{m}^{'})^{\sum\limits_{j\in A_{m-1}\setminus\{j_{0}\}}\frac{1}{\theta_j}} \Bigr)^{\theta_{m}} \right]^{\frac{1}{\theta_{m}}}<<2^{-n\alpha\delta_{m-1}}  \eqno (60) 
$$
$$
\times\left[\sum\limits_{0\leq s_{m}< \frac{n}{\gamma_{m}^{'}}}2^{-s_{m}\gamma_{m}^{'}(\frac{\gamma_{m}}{\gamma_{m}^{'}}-\delta_{m-1})\alpha\theta_{m}}(s_{m}+1)^{\lambda_{m}\theta_{m}}(n-s_{m}\gamma_{m}^{'}+1)^{\theta_{m}(\sum\limits_{j\in A_{m-1}}\lambda_{j}+\sum\limits_{j\in A_{m-1}\setminus \{j_{0}\}}\frac{1}{\theta_j})}\right]^{\frac{1}{\theta_{m}}} 
$$
under the condition  
 $\min\{\sum\limits_{j\in A_{m-1}\setminus \{j'\}}\lambda_{j} + \sum\limits_{j\in A_{m-1}\setminus \{j_{0}\}} \frac{1}{\theta_{j}}, \lambda_{j^{'}}+\frac{1}{\theta_{j^{'}}}\} > 0$.

If 
 $\frac{\gamma_{m}}{\gamma_{m}^{'}}-\delta_{m-1} >0$,  then  
$\delta_{m}=\min\{\frac{\gamma_{j}}{\gamma_{j}^{'}} : j=1,\ldots, m\}=\delta_{m-1}$, $A_{m}=\{j : \frac{\gamma_{j}}{\gamma_{j}^{'}}=\delta, j=1,\ldots, m\}=A_{m-1}$, $j_{1}=\min\{j : j\in A_{m}\}=j_{0}$, $j^{'}=\max\{j\in A_{m-1}\}$, $\overline{a}_{m-1}=(a_{1},...,a_{m-1})$.
Therefore 
 $\sum\limits_{j\in A_{m}\setminus\{j_{1}\}}\frac{1}{\theta_j}=\sum\limits_{j\in A_{m-1}\setminus\{j_{0}\}}\frac{1}{\theta_j}$. 

Now, by using Lemma 2 from (60), we obtain
$$
\sigma_{1}(n)<<2^{-n\alpha\delta_{m}}(n+1)^{\sum\limits_{j\in A_{m}}\lambda_{j}+\sum\limits_{j\in A_{m}\setminus\{j_{1}\}}\frac{1}{\theta_j}},
\eqno (61) 
$$
for 
 $\lambda_{m}\in \mathbb{R}$ and $\min\{\sum\limits_{j\in A_{m-1}\setminus \{j'\}}\lambda_{j} + \sum\limits_{j\in A_{m-1}\setminus \{j_{0}\}} \frac{1}{\theta_{j}}, \lambda_{j^{'}}+\frac{1}{\theta_{j^{'}}}\} > 0$ 

Let 
 $\frac{\gamma_{m}}{\gamma_{m}^{'}}-\delta_{m-1} <0$. Then  
  $\frac{\gamma_{m}}{\gamma_{m}^{'}}< \frac{\gamma_{j}}{\gamma_{j}^{'}}$, $j=1,...,m-1$. 
Since 
$\delta_{m-1}-\frac{\gamma_{m}}{\gamma_{m}^{'}} > 0$, by using Lemma 2 we have 
$$
\sum\limits_{0\leq s_{m}< \frac{n}{\gamma_{m}^{'}}}2^{-s_{m}\gamma_{m}^{'}(\frac{\gamma_{m}}{\gamma_{m}^{'}}-\delta_{m-1})\alpha\theta_{m}}(s_{m}+1)^{\lambda_{m}\theta_{m}}(n-s_{m}\gamma_{m}^{'}+1)^{\theta_{m}(\sum\limits_{j\in A_{m-1}}\lambda_{j}+\sum\limits_{j\in A_{m-1}\setminus \{j_{0}\}}\frac{1}{\theta_j})}
$$
$$
<<2^{-n(\frac{\gamma_{m}}{\gamma_{m}^{'}}-\delta_{m-1})\alpha\theta_{m}}(n+1)^{\lambda_{m}\theta_{m}}  \eqno (62)  
$$
under the condition  
$$
\min\{\sum\limits_{j\in A_{m-1}\setminus \{j'\}}\lambda_{j} + \sum\limits_{j\in A_{m-1}\setminus \{j_{0}\}} \frac{1}{\theta_{j}}, \lambda_{j^{'}}+\frac{1}{\theta_{j^{'}}}\} > 0.
$$
Now, from inequalities (60) and (62), we obtain
$$
\sigma_{1}(n) <<2^{-n\alpha\delta_{m}}(n+1)^{\lambda_{m}} \eqno (63)  
$$
in the case 
$\frac{\gamma_{m}}{\gamma_{m}^{'}}-\delta_{m-1} <0$ 
under the condition  $\lambda_{m}\in \mathbb{R}$ and   
 $$
 \min\{\sum\limits_{j\in A_{m-1}\setminus \{j'\}}\lambda_{j} + \sum\limits_{j\in A_{m-1}\setminus \{j_{0}\}} \frac{1}{\theta_{j}}, \lambda_{j^{'}}+\frac{1}{\theta_{j^{'}}}\} > 0.
 $$ 
 
Further, it follows from estimates (59) and (63) that 
$$
I_{n} <<2^{-n\alpha\delta_{m}}(n+1)^{\lambda_{m}}  \eqno (64) 
$$
in the case  
$\frac{\gamma_{m}}{\gamma_{m}^{'}}-\delta_{m-1} <0$ under the condition 
$\lambda_{m}\in \mathbb{R}$ and
$$
\min\{\sum\limits_{j\in A_{m-1}\setminus \{j'\}}\lambda_{j} + \sum\limits_{j\in A_{m-1}\setminus \{j_{0}\}} \frac{1}{\theta_{j}}, \lambda_{j^{'}}+\frac{1}{\theta_{j^{'}}}\} > 0. \eqno (65) 
$$ 

Let $\frac{\gamma_{m}}{\gamma_{m}^{'}}-\delta_{m-1} <0$, then 
$A_{m-1}\subset A_{m}$ and $j_{1}=\min{j\in A_{m}}=\min{j\in A_{m-1}}=j_{0}$.
Therefore, by using Lemma 2, we obtain
$$
\sum\limits_{0\leq s_{m}< \frac{n}{\gamma_{m}^{'}}}2^{-s_{m}\gamma_{m}^{'}(\frac{\gamma_{m}}{\gamma_{m}^{'}}-\delta_{m-1})\alpha\theta_{m}}(s_{m}+1)^{\lambda_{m}\theta_{m}}(n-s_{m}\gamma_{m}^{'}+1)^{\theta_{m}(\sum\limits_{j\in A_{m-1}}\lambda_{j}+\sum\limits_{j\in A_{m-1}\setminus \{j_{0}\}}\frac{1}{\theta_j})}
$$
$$
=
\sum\limits_{0\leq s_{m}< \frac{n}{\gamma_{m}^{'}}}(s_{m}+1)^{\lambda_{m}\theta_{m}}(n-s_{m}\gamma_{m}^{'}+1)^{\theta_{m}(\sum\limits_{j\in A_{m-1}}\lambda_{j}+\sum\limits_{j\in A_{m-1}\setminus \{j_{0}\}}\frac{1}{\theta_j})}
$$
$$
<<(n+1)^{\theta_{m}(\sum\limits_{j\in A_{m-1}}\lambda_{j}+\sum\limits_{j\in A_{m-1}\setminus \{j_{0}\}}\frac{1}{\theta_j})+\lambda_{m}\theta_{m} +1} = C(n+1)^{\theta_{m}(\sum\limits_{j\in A_{m}}\lambda_{j}+\sum\limits_{j\in A_{m}\setminus \{j_{1}\}}\frac{1}{\theta_j})} \eqno (66) 
$$
in the case  
 $\frac{\gamma_{m}}{\gamma_{m}^{'}}-\delta_{m-1} = 0$ under the condition 
 $\lambda_{m}+\frac{1}{\theta_{m}} > 0$ and (65).
  
 Now,  from inequalities (60) and (66) it follows that
$$
\sigma_{1}(n)<<2^{-n\alpha\delta_{m}}(n+1)^{\theta_{m}(\sum\limits_{j\in A_{m}}\lambda_{j}+\sum\limits_{j\in A_{m}\setminus \{j_{1}\}}\frac{1}{\theta_j})}
$$
in the case 
 $\frac{\gamma_{m}}{\gamma_{m}^{'}}-\delta_{m-1} = 0$ under the condition 
 $$
 \min\{\sum\limits_{j\in A_{m}\setminus \{j'\}}\lambda_{j} + \sum\limits_{j\in A_{m}\setminus \{j_{1}\}} \frac{1}{\theta_{j}}, \lambda_{m}+\frac{1}{\theta_{m}} > 0.
 $$

Therefore, from inequality (59) we obtain
$$
I_{n} <<\Bigl[2^{-n\alpha\delta_{m}}(n+1)^{\theta_{m}(\sum\limits_{j\in A_{m}}\lambda_{j}+\sum\limits_{j\in A_{m}\setminus \{j_{1}\}}\frac{1}{\theta_j})} +2^{-n\alpha\frac{\gamma_{m}}{\gamma_{m}^{'}}}(n+1)^{\lambda_{m}}\Bigr] 
$$
$$
<<2^{-n\alpha\delta_{m}}(n+1)^{\theta_{m}(\sum\limits_{j\in A_{m}}\lambda_{j}+\sum\limits_{j\in A_{m}\setminus \{j_{1}\}}\frac{1}{\theta_j})}
$$
in the case 
 $\frac{\gamma_{m}}{\gamma_{m}^{'}}-\delta_{m-1} = 0$ under the condition 
 $$
 \min\{\sum\limits_{j\in A_{m}\setminus \{j'\}}\lambda_{j} + \sum\limits_{j\in A_{m}\setminus \{j_{1}\}} \frac{1}{\theta_{j}}, \lambda_{m}+\frac{1}{\theta_{m}} > 0.
 $$
Lemma 3 is proved.

\begin{lemma}\label{lem 4}
{\it Let   $\overline{\gamma} = (\gamma_{1}, \ldots, \gamma_{m})$, $\overline{\varepsilon}=(\varepsilon_{1}, \ldots, \varepsilon_{m})$ and  $0<  \gamma_{j}$, $0< \varepsilon_{j}\leq \infty$, $\lambda_{j}$, $j=1,\ldots, m$, $m\geq 2$ and $\alpha \in (0, \infty)$, then  
$$
J_{n}=\left\|\left\{2^{-\alpha\langle\overline{s}, \overline{\gamma}\rangle}
\prod_{j=1}^{m}
(s_{j} + 1)^{\lambda_{j}}
\right\}_{\bar{s}\in \kappa^{m}(n, \bar{\gamma})}\right\|_{
l_{\bar \theta}}>>2^{-n\alpha}n^{\sum\limits_{j=1}^{m}\lambda_{j} +\sum\limits_{j=2}^{m}\frac{1}{\varepsilon_j}}.
$$
}
\end{lemma}   

{\bf Proof.} Let $m=2$. Then, by the definition of the set
$\kappa^{2}(n, \bar{\gamma}^{'})$ and Lemma 2 we have that
$$
J_{n}:=\left\|\left\{2^{-\alpha\langle\overline{s}, \overline{\gamma}\rangle}
\prod_{j=1}^{2}(s_{j} + 1)^{\lambda_{j}}
\right\}_{\bar{s}\in \kappa^{2}(n, \bar{\gamma})}\right\|_{
l_{\bar \theta}} \geq
$$
$$
\Bigl\{\sum\limits_{s_{2}< \frac{n}{\gamma_{2}}}\Bigl[\sum\limits_{s_{1}= \frac{n-s_{2}\gamma_{2}}{\gamma_{1}}} 2^{-\alpha\langle\overline{s}, \overline{\gamma}\rangle \varepsilon_{1}}\prod_{j=1}^{2}(s_{j} + 1)^{\lambda_{j}\varepsilon_{1}} \Bigr]^{\frac{\varepsilon_{2}}{\varepsilon_{1}}}  \Bigr\}^{\frac{1}{\varepsilon_{2}}}
$$
$$
= 2^{-n\alpha}\Bigl\{\sum\limits_{s_{2}< \frac{n}{\gamma_{2}}}(s_{2} + 1)^{\lambda_{2}\varepsilon_{2}}(n - s_{2}\gamma_{2} + 1)^{\lambda_{1}\varepsilon_{2}} \Bigr\}^{\frac{1}{\varepsilon_{2}}}>> 2^{-n\alpha}(n+1)^{\lambda_{1}+\lambda_{2}+\frac{1}{\varepsilon_{2}}}.   
$$
Now suppose that the lemma is true for $m-1 \geq 2$, that is,
$$
\left\|\left\{2^{-\alpha\langle\overline{s}_{m-1}, \overline{\gamma}_{m-1}\rangle}
\prod_{j=1}^{m-1}
(s_{j}+1)^{\lambda_{j}}
\right\}_{\bar{s}_{m-1}\in \kappa^{m-1}(n, \bar{\gamma})}\right\|_{
l_{\bar \theta}(\mathbb{Z}^{m-1}_{+})}
$$
$$
>>2^{-n\alpha}(n+1)^{\sum\limits_{j=1}^{m-1}\lambda_{j}+\sum\limits_{j=2}^{m-1}\frac{1}{\varepsilon_{j}}},
\eqno (67)  
$$
Let us prove the assertion of the lemma for $ m $.
By the definition of the set
$\kappa^{2}(n, \bar{\gamma}^{'})$, according to assumption (67) and Lemma 2, we have  
$$
J_{n} = 
$$
$$
\Biggl\{\sum\limits_{s_{m}< \frac{n}{\gamma_{m}}}2^{-s_{m}\gamma_{m}\alpha\varepsilon_{m}}(s_{m} + 1)^{\lambda_{m}\varepsilon_{m}}\left\|\left\{2^{-\alpha\langle\overline{s}_{m-1}, \overline{\gamma}_{m-1}\rangle}
\prod_{j=1}^{m-1}
(s_{j}+1)^{\lambda_{j}}
\right\}_{\bar{s}_{m-1}\in \kappa^{m-1}(n, \bar{\gamma})}\right\|_{
l_{\bar \theta}(\mathbb{Z}^{m-1}_{+})}^{\varepsilon_{m}} \Biggr\}^{\frac{1}{\varepsilon_{m}}}
$$
$$
>>2^{-n\alpha} \Biggl\{\sum\limits_{s_{m}< \frac{n}{\gamma_{m}}} (s_{m} + 1)^{\lambda_{m}\varepsilon_{m}}(n-s_{m}\gamma_{m}+1)^{(\sum\limits_{j=1}^{m-1}\lambda_{j}+\sum\limits_{j=2}^{m-1}\frac{1}{\varepsilon_{j}})\varepsilon_{m}} 
\Biggr\}^{\frac{1}{\varepsilon_{m}}}
$$
$$
>> 2^{-n\alpha}(n+1)^{\sum\limits_{j=1}^{m}\lambda_{j}+\sum\limits_{j=2}^{m}\frac{1}{\varepsilon_{j}}}.
$$
The lemma is proved.

\begin{rem} Lemma 3 for $\lambda_{j}=0$, $\theta_{j}=\theta$, $j=1,...,m$ and $\gamma_{j}=\gamma_{j}^{'}=1$ for $j = 1, \ldots, \nu $ and $1<\gamma_{j}^{'}<\gamma_{j}$ for $j=\nu + 1, ..., m$ was proved earlier in \cite[Lemma 2]{21}.
In the case $\lambda_{j}=0$,  $1\leq \theta_{j} < \infty$, $j=1,...,m$
and $\gamma_{j}=\gamma_{j}^{'}=1$ for $j=1, \ldots, \nu$ and $1<\gamma_{j}^{'}<\gamma_{j}$ for $j=\nu + 1, ..., m $ Lemma 3 is given in \cite{13}, \cite{14}.
For $\lambda_{j}=0$, $j=1,...,m$ from Lemma 3 and Lemma 4 we obtain \ cite [Lemma 2 and Lemma 3] {17}.
 \end{rem}

\section{Main results}\label{sec 2}

We now present the main results of the article.

\begin{theorem}\label{th 1} 
Let $\bar{p}=(p_{1},\ldots , p_{m})$, $\bar{q}=(q_{1},\ldots , q_{m})$, $\bar{\alpha}=(\alpha_{1},\ldots , \alpha_{m})$,  $\bar{\beta}=(\beta_{1},\ldots , \beta_{m})$, $\bar{\tau}^{(1)}=(\tau_{1}^{(1)},\ldots , \tau_{m}^{(1)})$, $\bar{\tau}^{(2)}=(\tau_{1}^{(2)},\ldots , \tau_{m}^{(2)})$, $\bar{\gamma}=(\gamma_{1},\ldots , \gamma_{m})$, $\bar{\gamma}^{'}=(\gamma_{1}^{'},\ldots , \gamma_{m}^{'})$, $\bar{r}=(r_{1},\ldots , r_{m})$, $\bar{\theta}=(\theta_{1},\ldots, \theta_{m})$ and 
$0< \theta_{j}\le \infty$, $1<\tau_{j}^{(1)}, \tau_{j}^{(2)} <+\infty,$
$1< p_{j} <q_{j} <+\infty,$ $\alpha_{j}, \beta_{j}\in \mathbb{R}$, $r_{j}> \frac{1}{p_{j}} - \frac{1}{q_{j}}$, $\gamma_{j}=\frac{r_{j} + \frac{1}{q_{j}} - \frac{1}{p_{j}}}{r_{j_{0}} + \frac{1}{q_{j_{0}}} - \frac{1}{p_{j_{0}}}}$, $1\leq \gamma_{j}^{'}\leq \gamma_{j}$,  $j=1,...,m$ and $r_{j_{0}} + \frac{1}{q_{j_{0}}} - \frac{1}{p_{j_{0}}}=\min\{r_{j} + \frac{1}{q_{j}} - \frac{1}{p_{j}} : j=1,...,m\}$,  $A=\min\{j : \frac{\gamma_{j}}{\gamma_{j}^{'}}=1, \,\, j=1,...,m\}$, $j_{1}=\min\{j\in A\}$.

If  
$$
\min\{\sum\limits_{j\in A\setminus \{j^{'}\}}(\beta_{j} - \alpha_{j}) +\sum\limits_{j\in A\setminus \{j_{1}\}} \left(\frac{1}{\tau_{j}^{(2)}}-\frac{1}{\theta_{j}}\right), \,\, \beta_{j^{'}} - \alpha_{j^{'}}+\frac{1}{\tau_{j^{'}}^{(2)}} - \frac{1}{\theta_{j^{'}}} \} > 0,
$$
then   
 $$
E_{n}^{(\bar{\gamma}^{'})}\left(S^{\bar r}_{\bar{p}, \bar\alpha, \bar{\tau}^{(1)}
, \bar{\theta}}B\right)_{\bar{q}, \bar\beta, \bar{\tau}^{(2)}}
\asymp
2^{^{^{-n\Bigl(r_{j_{0}} +\frac{1}{q_{j_{0}}}
-\frac{1}{p_{j_{0}}}\Bigr)}}}
n^{\sum\limits_{j\in A}(\beta_{j} - \alpha_{j}) +\sum\limits_{j\in A\setminus \{j_{1}\}}(\frac{1}{\tau_{j}^{(2)}}-\frac{1}{\theta
_{j}})_{+}}, 
$$
where  $j^{'}=\max\{j\in A\}$ and $a_{+}=\max\{a, 0\}$.
\end{theorem}

{\bf Proof.} Since   
$\varphi_{j}(t) =
t^{\frac{1}{p_{j}}}(1 + |\log_{2} t|)^{\alpha_{j}}$ and $\psi_{j}(t) =
t^{\frac{1}{q_{j}}}(1 + |\log_{2} t|)^{\beta_{j}}$ for  $t\in (0, 1]$ and $r_{j}> \frac{1}{p_{j}} - \frac{1}{q_{j}}$, $j=1,\ldots, m$, then the series  
$$
\sum\limits_{s_{j}=0}^{\infty}\Bigl(\frac{\psi_{j}(2^{-s_{j}})}
{\varphi_{j}(2^{-s_{j}})}2^{-s_{j}r_{j}}\Bigr)^{\varepsilon_{j}} = \sum\limits_{s_{j}=0}^{\infty} \Bigl(2^{-s_{j}(r_{j}+ \frac{1}{q_{j}}-\frac{1}{p_{j}})}(s_{j}+1)^{\beta_{j}-\alpha_{j}}\Bigr)^{\varepsilon_{j}}
$$
converge for $j=1,\ldots, m$. 

Therefore, by Theorem 3 \cite{16}, with $\overline{\gamma}$ replaced by $\overline{\gamma}^{'}$, we have
$$
E_{n}^{(\bar{\gamma}^{'})}(S_{\bar{p}, \bar\alpha, \bar{\tau}^{(1)}, \bar\theta}^{\bar r}B)_{\bar{q}, \bar\beta, \bar{\tau}^{(2)}}\leq
C\left\|\left\{\prod_{j=1}^{m}2^{-s_{j}(r_{j}+ \frac{1}{q_{j}}-\frac{1}{p_{j}})}(s_{j}+1)^{\beta_{j}-\alpha_{j}}\right\}_{\bar{s}\in Y^{m}(\bar{\gamma}^{'}, n)}
\right\|_{\bar\varepsilon}, \eqno (68) 
 $$
where 
 $\bar{\varepsilon}=(\varepsilon_{1},...,\varepsilon_{m}),$ 
$\varepsilon_{j}=\tau_{j}\beta_{j}',$  $\frac{1}{\beta_{j}}+
\frac{1}{\beta_{j}'}=1,$  $\beta_{j}=\frac{\theta_{j}}{\tau_{j}^{(2)}}.$
 
 Now applying Lemma 3 with $\lambda_{j}=\beta_{j}-\alpha_{j}$ and $\theta_{j}=\varepsilon_{j}$, $j=1,\ldots, m$,  $\alpha=r_{j_{0}} + \frac{1}{q_{j_{0}}} - \frac{1}{p_{j_{0}}}$ from (68) we get
  $$
E_{n}^{(\bar{\gamma}^{'})}(S_{\bar{p}, \bar\alpha, \bar{\tau}^{(1)}, \bar\theta}^{\bar r}B)_{\bar{q}, \bar\beta, \bar{\tau}^{(2)}}\leq
C2^{-n(r_{j_{0}} + \frac{1}{q_{j_{0}}} - \frac{1}{p_{j_{0}}})}(n+1)^{\sum\limits_{j\in A}(\beta_{j}-\alpha_{j})+\sum\limits_{j\in A\setminus \{j_{1}\}}(\frac{1}{\tau_{j}^{(2)}} -\frac{1}{\theta_{j}})}
$$
in the case 
 $1< \tau_{j}^{(2)} < \theta_{j} \leq \infty$, $j=1,\ldots , m$.

Let $1 \leq \theta_{j} \leq  \tau_{j}^{(2)} < \infty$, $j=1,\ldots , m$. Then it follows from the second statement of Theorem 3 \cite{16} that
$$
E_{n}^{(\bar{\gamma}^{'})}(S_{\bar{p}, \bar\alpha, \bar{\tau}^{(1)}, \bar\theta}^{\bar r}B)_{\bar{q}, \bar\beta, \bar{\tau}^{(2)}} \leq C  2^{-n(r_{j_{0}} + \frac{1}{q_{j_{0}}} - \frac{1}{p_{j_{0}}})\delta}(n+1)^{\sum\limits_{j=1}^{m}(\beta_{j}-\alpha_{j})}.
$$
The upper bound is proved.   
  
We now prove the lower bound.  
Let 
 $A=\{j: \gamma_{j}^{'}=\gamma_{j}, \,\, j=1,...,m\}$, $j_{1}=\min\{j\in A\}$ and $B = \{j: \tau_{j}^{(2)} <\theta_{j},\,\, j=1,...,m\}=\{1,...,m\}$. Then   
  $A\cap B\cup\{j_{1}\}=A$. We put  
   $\overline{s}^{0}=(s_{1}^{0},...,s_{m}^{0})$, where $s_{j}^{0}=s_{j}$ for $j\in A$ and $s_{j}^{0}=0$ for $j\notin A$.  

Consider the function  
$$
f_{1, n}(\overline{x})= n^{-\sum\limits_{j\in A\setminus \{j_{1}\}}\frac{1}{\theta_{j}}}\sum\limits_{\langle\bar{s}^{0}, \bar\gamma \rangle =n}\prod_{j=1}^{m}2^{-s_{j}^{0}(r_{j}+ 1-\frac{1}{p_{j}})}(s_{j}^{0}+1)^{-\alpha_{j}}\sum_{\bar{k}\in\rho(\bar{s}^{0})}e^{i\langle\bar{k},\bar{x}\rangle}.
$$
Then, according to continuity, the function   
$f_{1, n}\in L_{\bar{p}, \bar\alpha, \bar{\tau}^{(1)}}^{*}(\mathbb{T}^{m})$. 

The relation is known (see \cite{19}, \cite{20} ) 
$$
\Bigl
\|\sum_{\bar{k}\in\rho(\bar s)}e^{i\langle\bar{k},\bar{x}\rangle}
\Bigr\|_{\bar{p}, \bar\alpha, \bar{\tau}}^{*} = \prod_{j=1}^{m}\Bigl\|\sum_{k_{j}=2^{s_{j}-1}}^{2^{s_{j}}-1}e^{ik_{j}x_{j}}
\Bigr\|_{p_{j}, \alpha_{j}, \tau_{j}}^{*} \asymp \prod_{j=1}^{m}
2^{s_{j}(1-\frac{1}{p_{j}})}(s_{j}+1)^{\alpha_{j}}, \,\,       \eqno (69) 
 $$ 
for 
 $1<p_{j}, \tau_{j}<+\infty, \,\,$, $\alpha_{j}\in \mathbb{R}$, $s_{j}\in \mathbb{N}$, $ j=1,...,m$.

Taking into account relation (69), we obtain
$$
\Bigl
\|\Bigl\{\prod_{j=1}^{m}2^{s_{j}(r_{j}}\|\delta_{\bar s}(f_{1, n})\|_{\bar{p}, \bar\alpha, \bar{\tau}^{(1)}}^{*}\Bigr\}_{\bar{s}\in \mathbb{Z}_{+}^{m}}\Bigr\|_{l_{\bar\theta}} \asymp n^{-\sum\limits_{j\in A\setminus \{j_{1}\}}\frac{1}{\theta_{j}}}\Bigl
\|\Bigl\{\chi_{\kappa(\bar{\gamma}^{'}, n)}(\bar{s}^{0})\Bigr\}_{\langle\bar{s}^{0}, \bar{\gamma}^{'} \rangle = n}\Bigr\|_{l_{\bar\theta}},   \eqno (70) 
 $$
where 
  $\chi_{\kappa(\bar{\gamma}^{'}, n)}$---characteristic function of a set 
   $\kappa(\bar{\gamma}^{'}, n)= \{\bar{s}=(s_{1},...,s_{m})\in \mathbb{Z}_{+}^{m}: \,\, \langle\bar{s}, \bar{\gamma}^{'} \rangle = n \}$.

Let 
   $\tilde{\bar{s}}=(s_{j_{1}},...,s_{j_{|A|}})$, $\tilde{\bar{\gamma}}^{'}=(\gamma_{j_{1}},...,\gamma_{j_{|A|}})$, where $j_{i}\in A$, $i=1,...,|A|$ and $|A|$-- the number of elements of the set $A$.
 Then 
  $\langle\bar{s}^{0}, \bar{\gamma}^{'} \rangle = \sum\limits_{j \in A}s_{j}\gamma_{j}^{'} = \sum\limits_{i=1}^{|A|}s_{j_{i}}\gamma_{j_{i}}=\langle\tilde{\bar{s}}, \tilde{\bar{\gamma}}^{'}\rangle$.

Therefore, by Lemma 3 [17], relation (70) can be rewritten in the following form  
$$
\Bigl
\|\Bigl\{\prod_{j=1}^{m}2^{s_{j}r_{j}}\|\delta_{\bar s}(f_{1, n})\|_{\bar{p}, \bar\alpha, \bar{\tau}^{(1)}}^{*}\Bigr\}_{\bar{s}\in \mathbb{Z}_{+}^{m}}\Bigr\|_{l_{\bar\theta}} \asymp n^{-\sum\limits_{j\in A\setminus \{j_{1}\}}\frac{1}{\theta_{j}}}\Bigl
\|\Bigl\{\chi_{\kappa(\tilde{\bar{\gamma}}^{'}, |A|)}(\tilde{\bar{s}})\Bigr\}_{\langle\tilde{\bar{s}}, \tilde{\bar{\gamma}} \rangle = n}\Bigr\|_{l_{\bar\theta}} \asymp C_{1},   
 $$
where 
 $\kappa (\tilde{\bar{\gamma}}, |A|)=\{\tilde{\bar{s}}=(s_{j_{1}},...,s_{j_{|A|}}): \,\, \langle\tilde{\bar{s}}, \tilde{\bar{\gamma}} \rangle =n \}$ and 
   $\chi_{\kappa(\tilde{\bar{\gamma}}, |A|)}(\tilde{\bar{s}})$ ---characteristic function of a set 
    $\kappa(\tilde{\bar{\gamma}}, |A|)$.

Thus, the function  
 $F_{1, n}=C_{1}^{-1}f_{1, n} \in S_{\bar{p}, \bar\alpha, \bar{\tau}^{(1)}, \bar\theta}^{\bar r}B$.

By the definition of the best approximation of a function and by Theorem 2, we have 
$$
E_{n}^{(\bar{\gamma}^{'})}(F_{1, n})_{\bar{q}, \bar\beta, \bar{\tau}^{(2)}} \geq \|F_{1, n}\|_{\bar{q}, \bar\beta, \bar{\tau}^{(2)}} \geq C
\Bigl
\|\Bigl\{\prod_{j=1}^{m}2^{s_{j}(\frac{1}{\tau_{j}} - \frac{1}{q_{j}})}\|\delta_{\bar{s}^{0}}(f_{1, n})\|_{\bar\tau,  \bar{\tau}^{(2)}}^{*}\Bigr\}_{\langle\bar{s}^{0}, \bar{\gamma}^{'} \rangle = n}\Bigr\|_{l_{\bar{\tau}^{(2)}}}, 
$$
where  
$\bar\tau = (\tau_{1},...,\tau_{m})$ and $1< q_{j}< \tau_{j} < \infty$, $j=1,...,m$.
Hence, taking into account relation (69) and Lemma 3, we obtain 
$$
E_{n}^{(\bar{\gamma}^{'})}(F_{1, n})_{\bar{q}, \bar\beta, \bar{\tau}^{(2)}} >>
n^{-\sum\limits_{j\in A\setminus \{j_{1}\}}\frac{1}{\theta_{j}}}
\Bigl
\|\Bigl\{2^{-\langle\bar{s}^{0}, \bar{\gamma}^{'} \rangle(r_{j_{0}}+\frac{1}{q_{j_{0}} }- \frac{1}{p_{j_{0}}})}\Bigr\}_{\langle\bar{s}^{0}, \bar{\gamma}^{'} \rangle = n}\prod_{j=1}^{m}(s_{j}^{0} +1)^{\beta_{j} - \alpha_{j}}\Bigr\|_{l_{\bar{\tau}^{(2)}}} 
$$
$$
>>
n^{-\sum\limits_{j\in A\setminus \{j_{1}\}}\frac{1}{\theta_{j}}}
\Bigl
\|\Bigl\{2^{-\langle\tilde\bar{s}, \tilde\bar{\gamma} \rangle(r_{j_{0}}+\frac{1}{q_{j_{0}} }- \frac{1}{p_{j_{0}}})}\prod_{i=1}^{|A|}(\tilde{s_{j_{i}}} +1)^{\beta_{j} - \alpha_{j}}\Bigr\}_{\tilde\bar{s}\in \kappa(\tilde{\bar{\gamma}}^{'}, |A|)}\Bigr\|_{l_{\bar{\tau}^{(2)}}}.
$$
Now, using Lemma 3 for $\lambda_{j_{i}}=\beta_{j_{i}} - \alpha_{j_{i}}$, $i=1,...,|A|$, from this we obtain
$$
E_{n}^{(\bar{\gamma}^{'})}(F_{1, n})_{\bar{q}, \bar\beta, \bar{\tau}^{(2)}}>> 2^{-n(r_{j_{0}}+\frac{1}{q_{j_{0}} }- \frac{1}{p_{j_{0}}})}
(n +1)^{\sum\limits_{j\in A}(\beta_{j} - \alpha_{j})+\sum\limits_{j\in A\setminus \{j_{1}\}}(\frac{1}{\tau_{j}^{(2)}}-\frac{1}{\theta_{j}})}
$$
in the case  $1< \tau_{j}^{(2)} < \theta_{j} \leq \infty$, $j=1,\ldots , m$.

Let 
 $1 \leq \theta_{j} \leq  \tau_{j}^{(2)} < \infty$, $j=1,\ldots , m$. 
We will consider the function
$$
f_{2, n}(\overline{x})= \prod_{j=1}^{m}2^{-s_{j}(r_{j}+ 1-\frac{1}{p_{j}})}(s_{j}+1)^{-\alpha_{j}}\sum_{\bar{k}\in\rho(\bar{s})}e^{i\langle\bar{k},\bar{x}\rangle}
$$
for 
 $\langle\bar{s}, \bar\gamma^{'} \rangle \geq n$.
Then, using relation (69), we obtain that the function 
$F_{2, n}=C_{1}^{-1}f_{2, n} \in S_{\bar{p}, \bar\alpha, \bar{\tau}^{(1)}, \bar\theta}^{\bar r}B$ and 
$$
E_{n}^{(\bar{\gamma}^{'})}(F_{2, n})_{\bar{q}, \bar\beta, \bar{\tau}^{(2)}} \geq \|F_{2, n}\|_{\bar{q}, \bar\beta, \bar{\tau}^{(2)}} >>2^{-n(r_{j_{0}}+\frac{1}{q_{j_{0}} }- \frac{1}{p_{j_{0}}})}
(n +1)^{\sum\limits_{j\in A}(\beta_{j} - \alpha_{j})}.
$$
More generally, if the set 
$B = \{j: \tau_{j}^{(2)} <\theta_{j},\,\, j=1,...,m\}=\{1,...,m\}\neq \{j=1,...,m\}$, then we denote $B^{'}=A\cap B\cup\{j_{1}\}$. We put  
   $\overline{s}^{0}=(s_{1}^{0},...,s_{m}^{0})$, where $s_{j}^{0}=s_{j}$ for $j\in B^{'}$ and $s_{j}^{0}=0$ for $j\notin B^{'}$.  
Now, considering the function  
$$
f_{3, n}(\overline{x})= n^{-\sum\limits_{j\in A\setminus \{j_{1}\}}\frac{1}{\theta_{j}}}\sum\limits_{\langle\bar{s}^{0}, \bar\gamma \rangle =n}\prod_{j=1}^{m}2^{-s_{j}^{0}(r_{j}+ 1-\frac{1}{p_{j}})}(s_{j}^{0}+1)^{-\alpha_{j}}\sum_{\bar{k}\in\rho(\bar{s}^{0})}e^{i\langle\bar{k},\bar{x}\rangle}.
$$ 
 and arguing as in the previous cases, we can prove the lower bound.
The theorem is proved.

\begin{rem} In the case $\alpha_{j}=\beta_{j}=0$ and $p_{j}=\tau_{j}^{(1)}=p$,
$q_{j}=\tau_{j}^{(2)}=q$, $\theta_{j}=\theta$ for $j = 1,..., m$ from Theorem 1 follows the previously known results 
 by V.N. Temlyakov \cite[Theorem 2.2]{21}  and A.S. Romanyuk \cite[Theorem 2]{22}.  
For $\alpha_{j}=\beta_{j}=0$, $j=1,...,m$ and $\gamma_{j}^{'} =\gamma_{j}=1$ for $ j = 1, ..., \nu $ and $ \gamma_{j}^{'} <\gamma_{j}$, $ j = \nu + 1, ..., m$ Theorem 1 implies \cite[Theorem 2]{13} (also see \cite[Theorem 3.5]{14}), and for $\alpha_{j}=\beta_{j}=0$ and $\gamma_{j}^{'} \leq \gamma_{j}$ for $ j = 1, ...,m$ \cite[Theorem 1]{17}.
\end{rem}


This work was supported by a grant
  Ministry of Education and Science of the Republic of Kazakhstan (Project AP 08855579).


 \end{document}